\documentclass[11pt]{article}

\usepackage[a4paper,left=2cm,right=2cm,top=1.5cm,bottom=2cm]{geometry}
\usepackage[utf8]{inputenc}
\usepackage[T1]{fontenc}
\usepackage[english]{babel}
\usepackage{amsthm}
\usepackage{amsmath}
\usepackage{amssymb}
\usepackage{mathrsfs}
\usepackage{setspace}
\usepackage{float}
\usepackage{multirow}
\usepackage{stmaryrd}
\usepackage{color}
\usepackage{wrapfig}
\usepackage{pifont}
\usepackage{array}
\usepackage[colorlinks=true,linkcolor=ocre,citecolor=ocre]{hyperref}
\usepackage{dsfont}
\usepackage{upgreek}
\usepackage{nicefrac}
\usepackage{tikz}{}
\usepackage{gensymb}
\usepackage{FiraSans}
\usepackage{graphicx}
\usepackage{caption}
\renewenvironment{proof}{\noindent{\sffamily{\textbf{Proof :}}}}{\begin{flushright}$\square$\end{flushright}}

\newcommand{\IE}{\mathbb{E}}
\newcommand{\IN}{\mathbb{N}}

\newcommand{\IR}{\mathbb{R}}
\newcommand{\IC}{\mathbb{C}}
\newcommand{\IT}{\mathbb{T}}

\newcommand{\drm}{\mathrm d}

\newcommand{\CD}{\mathcal D}
\newcommand{\CS}{\mathcal S}

\newcommand{\CP}{\mathcal P}

\newcommand{\CR}{\mathcal R}
\newcommand{\CM}{\mathcal M}
\newcommand{\CQ}{\mathcal Q}
\newcommand{\CC}{\mathcal C}
\newcommand{\CH}{\mathcal H}

\newcommand{\CB}{\mathcal B}

\newcommand{\SL}{\mathscr{L}}

\newcommand{\IDC}{\mathds{1}}

\renewcommand{\P}{\mathsf{P}}

\newcommand{\PB}{\overline{\P}}
\newcommand{\PI}{\mathsf{\Pi}}

\newcommand{\DC}{\mathsf{C}}

\newcommand{\eps}{\varepsilon}

\newcommand{\Wick}[1]{\mathbf{:}#1\mathbf{:}}

\newcommand\bbC{\mathbb{C}}
\newcommand\bbE{\mathbb{E}}

\newcommand\bbN{\mathbb{N}}
\newcommand\bbR{\mathbb{R}}
\newcommand\bbP{\mathbb{P}}

\newcommand\bbZ{\mathbb{Z}}

\newcommand{\mcB}{\mathcal{B}}

\newcommand{\mcN}{\mathcal{N}}

\newcommand\mcS{\mathcal{S}}

\newcommand{\mcR}{\mathcal{R}}

\newcommand{\defeq}{:=}
\newcommand{\eqdef}{:=}

\definecolor{ocre}{RGB}{64,123,121}

\newcounter{item}
\numberwithin{item}{section}

\newtheorem{theorem}[item]{\sffamily Theorem}
\newtheorem{definition}[item]{\sffamily Definition}
\newtheorem{proposition}[item]{\sffamily Proposition}
\newtheorem{lemma}[item]{\sffamily Lemma}
\newtheorem{corollary}[item]{\sffamily Corollary}
\newtheorem{remark}[item]{\sffamily Remark}

\newtheorem*{theorem*}{\sffamily Theorem}
\newtheorem*{definition*}{\sffamily Definition}
\newtheorem*{proposition*}{\sffamily Proposition}
\newtheorem*{lemma*}{\sffamily Lemma}
\newtheorem*{corollary*}{\sffamily Corollary}
\newtheorem*{remark*}{\sffamily Remark}

\usepackage[explicit]{titlesec}
\titleformat{\section}{\centering\Large\bfseries}{\thesection \ --}{0.4em}{\Large\bfseries #1}
\titleformat{\subsection}{\centering\large\bfseries}{\thesubsection \ --}{0.4em}{\large\bfseries #1}
\titleformat{\subsubsection}{\centering\bfseries}{\thesubsubsection \ --}{0.4em}{\bfseries #1}

\let\emph\relax
\DeclareTextFontCommand{\emph}{\bfseries\em}

\usepackage{mathtools}
\numberwithin{equation}{section}
\mathtoolsset{showonlyrefs}

\begin{document}

\title{\bfseries Analysis of the Anderson operator}
\author{Ismaël BAILLEUL, Nguyen Viet DANG and Antoine MOUZARD}
\date{}

\maketitle
\abstract{We consider the continuous Anderson operator $H=-\Delta+\xi$ on a closed Riemannian compact surface $\mcS$. We provide a short self-contained functional analytic construction of the operator as an unbounded operator on $L^2(\mcS)$ based on its resolvent as a meromorphic family of operators. Our main result is a precise description of the Anderson heat semigroup $(e^{-tH})_{t>0}$ with two-sided Gaussian bounds and sharp Gaussian small time asymptotics for its kernel with a number of consequences on the spectrum of $H$.

Using these results, we introduce and study the associated Gaussian field that we call the Anderson Gaussian free field and prove that the law of its random partition function characterizes the law of the spectrum of $H$. We also give a construction of two measures on path space associated to the Anderson operator, the polymer measure and the ground state diffusion, as path in the random environment given by $\xi$. We relate the Wick square of the Anderson Gaussian free field to the renormalized occupation measure of a Poisson process of loops of diffusion paths and we further prove some large deviation results for the Anderson diffusion and its bridges.}
\vspace{0.5cm}

\section{Introduction}
\label{SectionIntro}

The Anderson operator is the Schrödinger operator
\begin{equation}
H=-\Delta+\xi
\end{equation}
with $\xi$ a Gaussian white noise, that is a random centered field with formal covariance given by
\begin{equation}
\IE\big[\xi(x)\xi(y)\big]=\delta_0(x-y).
\end{equation} 
In this work, we construct and study this operator on a closed Riemannian compact surface $\CS$ with metric $g$, measure $\mu$ and the Laplace-Beltrami  operator $\Delta$ as a {\it negative} operator ($\Delta=\partial_{x^1}^2+\partial_{x^2}^2$ on the flat $2$--torus $\mathbb{T}^2$). In this context, the white noise is an isometry from $L^2(\CS,\mu)$ to the space of random variable with finite variance such that
\begin{equation}
\IE\big[\langle\xi,\varphi\rangle\langle\xi,\psi\rangle\big]=\int_\CS\varphi(x)\psi(x)\mu(\drm x)
\end{equation}
for all $\varphi,\psi\in L^2(\CS)$. The difficulty lies in the roughness of the potential as $\xi\in C^{-1-\kappa}(\CS)$ for any $\kappa>0$, with $C	^{-1-\kappa}$ standing for the negative Besov-Hölder space of negative regularity $-1-\kappa$. Fukushima and Nakao \cite{FN77} constructed the operator $H$ on a segment $[0,L]$ with $\xi=\drm B$ the distributional derivative of a Brownian motion, via Dirichlet form methods. In higher dimension the study of this operator falls in the range of singular SPDEs which received a lot of attention over the last decade with the introduction of regularity structures by Hairer \cite{Hai14} and paracontrolled calculus by Gubinelli, Imkeller and Perkowski \cite{GIP}. The first construction of the Anderson operator in this singular framework was achieved by Allez and Chouk \cite{AllezChouk} on the torus $\IT^2$ using paracontrolled calculus. It was then generalized to different contexts in \cite{GUZ,Labbe,Mouzard} including boxes $[0,L]^d$ for $d\in\{2,3\}$ with different boundary conditions and compact surfaces using both regularity structures and paracontrolled calculus. See also \cite{MvZ,MouzardOuhabaz} for an approach based on the Dirichlet form using an exponential transform. 

\medskip

The main difficulty for the construction of the Anderson operator is that for any smooth function $u$ one has
\begin{equation}
Hu=-\Delta u+u\xi\in C^{-1-\kappa}(\CS)
\end{equation}
where $\xi$ is irregular everywhere on the surface $\CS$. In order to gain some regularity for $Hu$, the idea is to induce roughness in the function $u$ depending on $\xi$, in such a way that $-\Delta u$ cancels out the roughest part of the product $u\xi$. This is precisely where the singularity appears in two dimensions: the natural regularity for these functions is $u\in C^{1-\kappa}$, hence the singularity of the product $u\xi$. This is where regularity structures or paracontrolled calculus appear, one considers a class of functions $u$ with some prescribed local behavior depending on the noise which allows to gain regularity in $Hu$ up to a probabilistic renormalization procedure. This path was followed by the recent constructions in \cite{AllezChouk,GUZ,Labbe,Mouzard} and the operator is constructed as a limit
\begin{equation}
H=\lim_{\eps\to0}\big(-\Delta+\xi_\eps-c_\eps\big)
\end{equation}
with $\xi_\eps$ a regularization of the noise and $c_\eps$ an explicit divergent quantity, logarithmic in $\varepsilon$ in two dimensions. In this work, we propose a construction in the spirit of the initial work of Allez and Chouk \cite{AllezChouk} with the additional input of the meromorphic Fredholm Theorem. We construct a family of operators $R_H(z)$ as a meromorphic family of operators that corresponds to the resolvent
\begin{equation}
R_H(z) = (H-z)^{-1}
\end{equation}
with poles the discrete eigenvalues of $H$.

 We summarize our results on the construction of $H$ and its basic properties as the following light Theorem: 
\begin{theorem*}
We construct the Anderson operator as a limit
\begin{align*}
H:=\lim_{\varepsilon\rightarrow 0^+}\left(-\Delta+\xi_\varepsilon-\frac{\vert \log(\varepsilon)\vert}{4\pi} \right),
\end{align*}
where $H$ is unbounded self--adjoint operator acting on $L^2(\mathcal{S})$,
with compact resolvent $(H-z)^{-1}=\lim_{\varepsilon\rightarrow 0^+}\left(-\Delta+\xi_\varepsilon-\frac{\vert \log(\varepsilon)\vert}{4\pi}-z \right)^{-1}$. The operator $H$ has therefore a 
\textbf{discrete spectrum} $(\lambda_n)_{n\in \mathbb{N}}$ with a corresponding sequence of $L^2$ normalized eigenfunctions $(v_n)_{n\in \mathbb{N}}$.
\end{theorem*}

 We then provide a precise description of the Anderson heat semigroup $(e^{-tH})_{t>0}$ including two-sided Gaussian bounds, Schauder estimates as well as small time asymptotic comparison with respect to the usual heat kernel $(e^{t\Delta})_{t>0}$. Most results of the present article follow from our main Theorem \ref{ThmAsymptotic}.
Let us state a baby version of Theorem \ref{ThmAsymptotic} and its important Corollary \ref{CorAsymptoticKernel} to give the reader a flavour of our results~:

\medskip

\begin{theorem*}
For $\alpha\in[0,2)$ and $\beta\in[-1,1]$, the Anderson heat operator $e^{-tH}$ satisfies Schauder estimates of the form 
\begin{equation}
\big\| e^{-tH}u_0 \big\|_{C^\alpha}\lesssim  t^{-\frac{\alpha-\beta}{2}}  \|u_0\|_{C^\beta}
\end{equation}
for $\beta<\alpha$. For $\alpha\in[0,1)$ and $\beta\in(-1,1]$, we compare the Anderson heat flow $e^{-tH}$ with the heat operator $e^{t\Delta}$, we get some \textbf{comparison estimates} on $e^{-tH}-e^{t\Delta}$ of the form
\begin{equation}
\big\| e^{-tH}u_0 - e^{t\Delta}u_0 \big\|_{C^\alpha}\lesssim t^{ -\frac{\alpha-\beta}{2}+\rho}\|u_0\|_{C^\beta}
\end{equation}
for any $\rho<\frac{1}{2}$.
\end{theorem*}

\medskip
Let us make some quick comments on the numerology of the above result. Schauder estimates give an explosion of order $t^{-\frac{\alpha-\beta}{2}}$ as $t$ goes to $0$ for both $\|e^{-tH}(u_0)\|_{C^\alpha} $ and $\|e^{t\Delta}(u_0)\|_{C^\alpha}$, but for the difference $e^{-tH}-e^{t\Delta}$ of the semigroups, one \textbf{gains a factor} $t^{\rho}$ in the blow-up of $\|e^{-tH}(u_0) - e^{t\Delta}(u_0)\|_{C^\alpha}$ when $t>0$ goes to $0$.

We give a number of consequences of Theorem \ref{ThmAsymptotic} for the spectrum of the Anderson operator, this recovers in particular the Weyl law for the Anderson operator
\begin{equation}
\big|\big\{\lambda\in\sigma(H)\,;\,\lambda\leq a\big\}\big| \underset{a\to+\infty}{\sim} \frac{\mu(\mcS)}{4\pi}\,a
\end{equation}
obtained previously in \cite{Mouzard}.

\medskip

We then introduce and study the Anderson Gaussian Free Field $\phi$, that is the centered Gaussian field with covariance given by the Green function of $H$. Our Gaussian upper bound on the Anderson heat kernel implies that its Green function has a logarithmic divergence along the diagonal, as for the Laplace-Beltrami operator on $\CS$, hence this field falls in the range of log-correlated fields. However, the particularity here is that we work with two layers of randomness as this is a Gaussian random field given the random environment $\xi$ on $\CS$. Like the usual Gaussian free field, it takes values in $C^{-\delta}(\CS)$ for any $\delta>0$. We construct its Wick square $\Wick{\phi^2}$ via a renormalization procedure, almost surely given our random singular environment of the Anderson operator. We relate its random partition function
\begin{equation}
Z(\lambda)=\IE\Big[e^{-\lambda\langle\Wick{\phi^2},1\rangle}\Big]
\end{equation}
to the determinant of the Anderson resolvent, this is the content of Theorem \ref{ThmWickCharacterizationSpectrumH}. In particular, this gives a characterization of the law of the spectrum of $H$. This field is of independent interest in the study of random Gibbs measure for stochastic PDEs, see for example the recent works \cite{BVZ,EMR}.

\medskip

Finally we consider two natural path measures associated to the Anderson operator, that is the polymer measure and the diffusion associated to the Anderson ground state. Our precise description of the Anderson heat kernel yields a number of properties of these measures similar to the Wiener measure. This has to be considered in perspective with the fact that the random polymer measure is almost surely singular with respect to the Wiener measure. This fact was already proved by Cannizzaro and Chouk \cite{CannizzaroChouk} where they constructed the measure via a KPZ equation and a Girsanov transform. We follow here a more direct approach since the Anderson heat kernel gives the probability transition of the underlying path, in the spirit of the construction in \cite{AKQ} of the random continuum polymer by Albert, Khanin and Quastel with a space-time white noise in dimension $1+1$. To conclude our paper, we relate the renormalized Wick square of the Anderson Gaussian free field with the renormalized occupation measure of a Poisson process of loops of Anderson diffusion paths in the spirit of the work of Le Jan.

\medskip

In Section \ref{SectionTools}, we explain the tools needed for our construction and study of the Anderson operator. We give the definition and needed basic properties of the paracontrolled calculus as well as the meromorphic Fredholm Theorem. In Section \ref{SectionResolventAnderson}, we construct the meromorphic family $R_H(z)$ via a renormalization procedure using the paracontrolled calculus. In Section \ref{SectionHeatKernelAsymptotics}, we study the Anderson heat semigroup $(e^{-tH})_{t>0}$. This is where our main theorems and consequences for the Anderson operator are stated and proved. In Section \ref{SectionAndersonGFF}, we introduce and study the Anderson Gaussian free field while the associated polymer and the Anderson diffusion are studied in Section \ref{SectionPolymer}. Finally, we respectively give in Appendix \ref{SectionFredholm} and \ref{SubsectinoGeometricLP} context on the Meromorphic Fredholm theory with a parameter and the Geometric Littlewood-Paley decomposition.

\bigskip

\noindent {\bf Acknowledgments --} We would like to thank the anonymous reviewer for her/his patience and for pointing out a mistake in the previous proof of our main theorem as well as Hugo Eulry, Tristan Robert and Immanuel Zachhuber for useful discussions on several points of the present work. N.V.D. would like to thank the Institut Universitaire de France for support.

\section{Tools for the analysis}
\label{SectionTools}

In this section, we recall what we need from paraproducts and meromorphic Fredholm theory. For a complete description of paracontrolled calculus we refer the reader to \cite{BCD, GIP, BB2, Mouzard}; we will only use here what we recall below.

\medskip

$\rhd$ {\it Paraproduct and co.} For any distribution $f$ on the $d$-dimensional torus $\IT^d$ one can consider the Littlewood-Paley projector
\begin{equation*}
\big(\Delta_nf\big)(x):=2^{d(n-1)}\int_{\IR^d}\chi\big(2^{n-1}(x-y)\big) f(y) \drm y
\end{equation*} 
with $\chi\in\CS(\IR^d)$ and $\text{supp}\ \widehat\chi\subset\{\frac{1}{2}\le|z|\le 2\}$ for $n\ge1$ approximately localized in frequency space in annuli of size $2^n$, and
\begin{equation*}
\big(\Delta_0f\big)(x):=\int_{\IR^d}\chi_0(x-y) f(y) \drm y
\end{equation*}
with $\chi_0\in\CS(\IR^d)$ and $\text{supp}\ \widehat\chi_0\subset\{|z|\le 1\}$. This allows to decompose $f$ as a sum of smooth functions
$$
f = \sum_{n\geq0} \Delta_nf
$$
and to measure the regularity with the Besov spaces associated to the norm
$$
\|f\|_{B_{p,q}^\alpha} := \bigg(\sum_{n\ge0}2^{\alpha pn}\|\Delta_nf\|_{L^q(\IT^d)}^p\bigg)^{\frac{1}{p}}
$$
for $p,q\in[1,+\infty]$ and $\alpha\in\IR$. One recovers the usual Hölder spaces for $p=q=\infty$ and $\alpha\in\IR^+\backslash\IN$ and the Sobolev spaces for $p=q=2$ and $\alpha\in\IR\backslash\IN$. In the following we consider the Besov-Hölder spaces $C^\alpha=B_{\infty,\infty}^\alpha$ and the Besov-Sobolev spaces $H^\alpha=B_{2,2}^\alpha$. This also allows to decompose formally the product of two distributions into
\begin{align*}
fg &= \sum_{n<m-1}\Delta_nf\Delta_mg + \sum_{m<n-1}\Delta_nf\Delta_mg + \sum_{\vert n-m\vert\leq 1}\Delta_nf\Delta_mg\\
&=: P_fg+P_gf+\Pi(f,g)
\end{align*}
with the paraproducts $P_fg$ and $P_gf$ being always well--defined on the spaces of distributions on the torus, the possible singularity being encoded in the resonant term $\Pi(f,g)$. While this paraproduct $P$ was introduced by Bony \cite{Bony}, this can be adapted in our $2$-dimensional manifold setting where one can decompose the product of any two smooth functions $f,g$ on $\mcS$ as
\begin{equation} \label{EqParaproductResonantIdentity}
fg = {\sf P}_fg + {\sf P}_gf + {\sf \Pi}(f,g)
\end{equation}
using some paraproduct and resonant operators $\sf P$ and $\sf \Pi$ with the following continuity properties. See Appendix \ref{SubsectinoGeometricLP} for the definition of the paraproduct and resonant product.

\begin{proposition}
\begin{enumerate}
   \item[(i)] For any $\alpha_1, \alpha_2\in\bbR$, the paraproduct operator 
   $$
   {\sf P} : (f,g)\mapsto {\sf P}_fg,
   $$ 
   maps continuously $C^{\alpha_1}(\mcS)\times C^{\alpha_2}(\mcS)$ into $C^{\alpha_1\wedge 0 + \alpha_2}(\mcS)$. For $\alpha_1\neq0$, it also maps continuously the space $C^{\alpha_1}(\mcS)\times H^{\alpha_2}(\mcS)$ and $H^{\alpha_1}(\mcS)\times C^{\alpha_2}(\mcS)$ into $H^{\alpha_1\wedge 0 + \alpha_2}(\mcS)$.

   \item[(ii)] For any $\alpha_1, \alpha_2\in\bbR$ such that $\alpha_1+\alpha_2>0$, the resonant operator 
   $$
   {\sf \Pi} : (f,g)\mapsto{\sf \Pi}(f,g),
   $$ 
   maps continuously $C^{\alpha_1}(\mcS)\times C^{\alpha_2}(\mcS)$ into $C^{\alpha_1+\alpha_2}(\mcS)$, and it sends continuously $C^{\alpha_1}(\mcS)\times H^{\alpha_2}(\mcS)$ into $H^{\alpha_1+\alpha_2}(\mcS)$.  
\end{enumerate}
\end{proposition}

\medskip

The decomposition \eqref{EqParaproductResonantIdentity} thus makes sense for all $f\in C^{\alpha_1}(\mcS), g\in C^{\alpha_2}(\mcS)$ provided $\alpha_1+\alpha_2>0$ -- this is usually refered to as Young's condition, which ensures that the product of two distributions is well-defined (it is a necessary condition). The reader will find more details on these paraproduct and resonant operators in Appendix \ref{SubsectinoGeometricLP}. In the context of singular SPDEs, these tools were used in the seminal work by Gubinelli, Imkeller and Perkowski \cite{GIP} where they consider the corrector
\begin{equation*} 
{\sf C}(a,b,c) \defeq {\sf \Pi}({\sf P}_ab,c) - a{\sf \Pi}(b,c)
\end{equation*}
for smooth functions $a,b,c$. Its extension to a manifold setting was worked out in Bailleul and Bernicot's work \cite{BB1} in a general parabolic setting and in Mouzard's work \cite{Mouzard} for the mixed elliptic Sobolev and Hölder regularities.

\begin{proposition} \label{PropCorrector}
For any $\alpha\in(0,1)$ and $\alpha_2,\alpha_3\in\IR$ such that $\alpha_1+\alpha_2+\alpha_3>0>\alpha_2+\alpha_3$, the trilinear operator $\DC$ is continuous from $C^{\alpha_1}(\mcS)\times C^{\alpha_2}(\mcS)\times C^{\alpha_3}(\mcS)$ into $C^{\alpha_1+\alpha_2+\alpha_3}(\mcS)$ and from $H^{\alpha_1}(\mcS)\times C^{\alpha_2}(\mcS)\times C^{\alpha_3}(\mcS)$ into $H^{\alpha_1+\alpha_2+\alpha_3}(\mcS)$.
\end{proposition}

Since $\alpha_2+\alpha_3<0$, the resonant product $\PI(b,c)$ is singular but this is also the case of the product $\PI(\P_ab,c)$, since $\alpha_1>0$. The previous continuity estimates on the corrector states that the singular part of each terms cancel each other given that $a$ is regular enough, that is $\alpha_1+\alpha_2+\alpha_3>0$. This is the backbone of the resolution of a number of singular SPDEs within paracontrolled calculus.

\medskip

Given $z_0\notin\sigma(\Delta)$ we will use occasionally the paraproduct-like operator $\overline{\sf P}$ defined by the intertwining relation
$$
\overline{\sf P}_fg \defeq (-\Delta + z_0)^{-1}{\sf P}_f(-\Delta + z_0)g,
$$
following some ideas from \cite{BB2} in the parabolic setting and \cite{Mouzard} in the elliptic setting. In particular, the operator $\PB$ has the same continuity properties as the operator $\P$. This operator $\overline{\sf P}$ depends on $z_0$, which will be fixed throughout, so we do not record it in the notation to lighten the redaction. In particular, it enjoys some continuity estimates whose constants are uniform for $z_0\ge 1$. It will also be important for us that the associated corrector
\begin{equation} \label{EqModifiedCorrector}
\overline{\sf C}(a,b,c) \defeq {\sf \Pi}(\overline{\sf P}_ab,c) - a{\sf \Pi}(b,c)
\end{equation}
enjoys the same continuity property as $\sf C$, stated in Proposition \ref{PropCorrector}, with $z_0$-uniform constants for $z_0\geq 1$. The space white noise $\xi$ is the centered Gaussian distribution with formal covariance given by a delta function
\begin{equation*}
\IE\big[\xi(x)\xi(y)\big]=\delta_x(y)
\end{equation*}
for any $x,y\in\CS$. One can then prove that it belongs almost surely to $C^{-1-\kappa}(\CS)$ for any $\kappa>0$ in two dimensions, so the product
\begin{equation*}
\xi f=\P_\xi f+\P_f\xi+\PI(\xi,f)\in C^{-1-\kappa}
\end{equation*}
is well-defined only for $f\in H^\alpha$ with $\alpha>1$. The estimates on paraproducts give in this case that the sum $\P_\xi f+\PI(\xi,f)$ is of better regularity than the product, $\P_f \xi$ being the roughest part. Motivated by the ideas from the paracontrolled calculus from \cite{GUZ,Mouzard}, we introduce the map
\begin{equation*}
\Phi(f):=f+(-\Delta+z_0)^{-1}(\P_f\xi)
\end{equation*}
which rewrites via the intertwining relation
\begin{equation*}
\Phi(f)=f+\PB_f X
\end{equation*}
with 
\[
X:=(-\Delta+z_0)^{-1}(\xi).
\]
The map $\Phi$ is a perturbation of the identify, thus invertible if the perturbation is small enough. This will be achieved by taking $z_0$ large enough, depending on $\xi$, with the following lemma. It follows from the paraproduct's continuity estimates and some estimates on the resolvent of the Laplace-Beltrami operator. For the rest of the work we fix a parameter $\kappa>0$ which can be taken arbitrarily small such that $\xi\in C^{-1-\kappa}$ hence $X\in C^{1-\kappa}$.

\medskip

\begin{lemma} \label{lem:adjustz0}
For any $\gamma\in \mathbb{R}$ and $\eta>0$ there exists $m=m(\eta,\kappa,\xi)>0$ such that for every real parameter $z_0\geq m$ one has
$$
\big\Vert(-\Delta + z_0)^{-1}\big\Vert_{\mcB(H^\gamma,H^{\gamma+2})} < \eta
$$
as well as
$$
\big\| f\mapsto (-\Delta + z_0)^{-1}(\P_f\xi) \big\|_{\mcB(H^\gamma,H^{\gamma\wedge 0 + 1-\kappa})}<1.
$$
\end{lemma}

\medskip

We get that $\Phi:H^\gamma\to H^\gamma$ is invertible for any $\gamma<1-\kappa$, we denote by $\Gamma$ its inverse. It is defined by the implicit equation
\begin{equation*}
\Gamma (g) = \PB_{\Gamma(g)}X + g,
\end{equation*}
which is a first order paracontrolled expansion. While our choice of parameter is different, the idea to take a truncation depending on the noise to ensure that $\Phi$ is invertible goes back to Gubinelli, Ugurcan and Zachhuber in \cite{GUZ}. This was later generalized to a second order paracontrolled $\Gamma$ map by the third author in \cite{Mouzard}, this is not needed for our work.

\medskip

$\rhd$ {\it Meromorphic Fredholm theory.} Finally, analytic Fredholm theory provides some conditions under which one can invert a family of Fredholm operators acting on some Hilbert space. Let $U$ be a connected open subset of the complex plane $\bbC$ and $(\mathcal{H},\langle\cdot,\cdot\rangle)$ be a Hilbert space. A family $\big(A(z)\big)_{z\in U}$ of linear maps from $\mathcal{H}$ into itself is said to be holomorphic if and only if the map $A$ is $\bbC$-differentiable in $U$. This is equivalent to requiring that the $\bbC$-valued function 
$$
z\mapsto \left\langle y , A(z)x\right\rangle
$$ 
is holomorphic for any $x,y$ in $\mathcal{H}$. The family $\big(A(z)\big)_{z\in U}$ is said to be {\it finitely meromorphic} if for any $z\in U$ there exists a finite collection of operators $(A_j)_{1\leq j\leq n_0}$ {\it of finite rank} and a holomorphic family $A_0(\cdot)$ such that one has
$$
A(z') = A_0(z') + (z'-z)^{-1}A_1 + \cdots + (z'-z)^{-n_0}A_{n_0}
$$
for $z'$ near $z$. In particular, this implies that the poles are isolated of finite order. We shall need a version with parameters of the meromorphic Fredholm Theorem where $A(z,e)$ depends continuously on a parameter $e$ taken in a metric space $(E,d)$.

\smallskip

\begin{theorem} \label{ThmFredholm}
Let $U\subset\bbC$ be a connected open subset of the complex plane. Let $(E,d)$ be a metric space and $\big(K(z,e)\big)_{z\in U, \, e\in E}$ a finitely meromorphic family of compact operators depending continuously on $e\in E$. If for every $v_0\in E$, the operator $\big(\textrm{Id} - K(z_0,e)\big)^{-1}$ exists at some point $z_0\in U$ for all $e$ in a neighborhood of $v_0$ then the family 
$$
(z'\in U)\mapsto \big(\textrm{Id} - K(z',e)\big)^{-1}
$$ 
is a well-defined meromorphic family of operators with poles of finite rank which depends continuously on $e\in E$.
\end{theorem}

\medskip

A proof of this statement is given in Appendix \ref{SectionFredholm}. To conclude this section we recall that a sequence $(h_n)_{n\geq 0}$ of Banach space-valued meromorphic functions defined on a common open subset of $\bbC$ converges to a limit meromorphic function $h$ if $h_n$ converges uniformly to $h$ on every compact set that does not contain any pole of $h$.

\section{Construction of the Anderson operator}
\label{SectionResolventAnderson}

In this section we construct the Anderson operator
$$
H=-\Delta+\xi
$$
where $\xi$ stands for a space white noise and $\Delta$ is the non-positive Laplace-Beltrami operator associated with the Riemannian metric on $\mcS$. One can construct $\xi$ as a random series 
$$
\sum_{n\geq 0} \gamma_n f_n
$$ 
where the $f_n$ are the eigenfunctions of $\Delta$ and the $\gamma_n$ are a family of independent centered Gaussian random variables with unit variance -- this goes back to Paley and Zygmund with Fourier series \cite{PaleyZygmund30,PaleyZygmund32}. The main difficulty of the construction of $H$ lies in the roughness of the noise, and its construction in two dimensions involves a renormalization procedure. 

\smallskip

We construct the unbounded operator $H$ on $L^2(\mcS)$ by its resolvent map $R_H(z)$ as a meromorphic function of $z$ with values in operators. A reader already familiar with one of the previous constructions of the Anderson operator \cite{AllezChouk, Labbe, GUZ, Mouzard} may skip this section and keep in mind that our point of view is to construct the resolvent of this operator as a meromorphic function on $\bbC$. This is in the spirit of the original construction by Allez and Chouk \cite{AllezChouk} where they consider the resolvent on $[z_0,+\infty)\subset\IC$ for $z_0=z_0(\xi)>0$ large enough. Given an operator $A$, denote its resolvent as
$$
R_A(z)=(A-z)^{-1}
$$
which is well-defined for $z\notin\sigma(A)$. For some operators $A$ and $B$ we have the identity
\begin{equation*}
R_A(z)-R_B(z) = R_A(z)(B-A)R_B(z),
\end{equation*}
which is often refered to as the second resolvent identity.  Applying this with the Laplace-Beltrami and the Anderson operators formally gives
\begin{equation*}
(-\Delta+z_0)^{-1}-R_H(z)=(-\Delta+z_0)^{-1}(\xi-z_0-z) R_H(z).
\end{equation*}

In the sequel, we will use the following shorthand notation. For any linear operator $A:C^\infty(\mathcal{S})\mapsto \mathcal{D}^\prime(\mathcal{S})$ we denote by
${\sf P}_A\xi$ and ${\sf P}_\xi A$ the linear operators $\varphi\in C^\infty(\mathcal{S})\mapsto {\sf P}_{A\varphi}\xi$ and $\varphi\in C^\infty(\mathcal{S})\mapsto {\sf P}_\xi(A\varphi) $. This intuitive notation is very helpful to make a number of computations look short.

\smallskip

Since the operator composition $\xi R_H(z)$ will be singular in our case, we decompose it as
\begin{equation*}
\xi R_H(z)=\P_{R_H(z)}\xi+\P_{\xi}R_H(z)+\PI\big(\xi,R_H(z)\big)
\end{equation*}
where the roughest term is given by the first paraproduct. We inject this decomposition in the resolvent identity, which yields
\begin{equation} \label{EqResolvantFixedPoint} \begin{split}
R_H(z)&=-\left(-\Delta+z_0 \right)^{-1}\P_{R_H(z)} \xi+R_H(z)^\sharp\\
&=-\PB_{R_H(z)}\left(-\Delta+z_0 \right)^{-1}\xi+R_H(z)^\sharp\\
&=-\PB_{R_H(z)}X+R_H(z)^\sharp,
\end{split} \end{equation}
using the $z_0$-dependent paraproduct $\PB$, where 
\begin{equation*}
R_H(z)^\sharp=  \left(-\Delta+z_0 \right)^{-1} \bigg(\text{Id} - {\sf P}_\xi R_H(z)-\PI\big(\xi,R_H(z) \big)+(z+z_0)R_H(z) \bigg) 
\end{equation*}
is of better formal regularity. Following the ideas of paracontrolled calculus, our goal is to construct the resolvent given a renormalized product. In relation with the previous construction of the Anderson operator, the Ansatz on $R_H(z)$ imposes that it takes values in a space of functions paracontrolled at first order by $X$. The term $R_H(z)^\sharp$ involves a resonant product for which we use the paracontrolled ansatz, that is
\begin{align*}
\PI\big(\xi,R_H(z)\big)&=  -  \PI\big(\xi,\PB_{R_H(z)}X\big) + \PI\big(\xi,R_H(z)^\sharp\big) \\
&=-\PI(\xi,X)R_H(z) -  \overline\DC\big(R_H(z),X, \xi\big)+\PI\big(\xi,R_H(z)^\sharp\big)
\end{align*}
with $X=\left(-\Delta+z_0 \right)^{-1}(\xi)\in C^{1-\kappa}(\CS)$ where the last identity involves the corrector $ \overline\DC$. The resonant term $\PI(\xi,X)$ is a priori ill-defined and has to be renormalized as a Wick product which gives an element of $C^{-2\kappa}(\CS)$: this is the content of Proposition \ref{PropCvgce} below. 

\smallskip

Back to \eqref{EqResolvantFixedPoint}, this is where the map $\Gamma$ appears naturally since
\begin{equation*}
R_H(z)=-\PB_{R_H(z)}X+R_H(z)^\sharp\quad\iff\quad R_H(z)=\Gamma\big(R_H(z)^\sharp\big)
\end{equation*}
for $z_0$ large enough depending on $\xi$. The resolvent of the Anderson operator is then given by the solution to the equation 
\[
R_H(z)=\Gamma\big(R_H(z)^\sharp\big)
\] 
with
\begin{align*}
&R_H(z)^\sharp =   \\
&\left(-\Delta+z_0 \right)^{-1}  \left(\text{Id}-P_\xi R_H(z)+\mcR\PI(\xi,X)R_H(z) +  \overline\DC\big(R_H(z),X, \xi\big)-\PI\big(\xi,R_H(z)^\sharp\big) +(z+z_0)R_H(z) \right) 
\end{align*}
with $\mcR\PI(\xi,X)$ given as an element of $C^{-2\kappa}(\CS)$. For its renormalization, consider the heat regularized space white noise
$$
\xi_\eps \defeq e^{\eps\Delta}(\xi)
$$
with $\eps>0$ as well as
$$
X_\eps \defeq (-\Delta + z_0)^{-1}(\xi_\eps).
$$ 
The resonant term $\PI(\xi_\eps,X_\eps)$ diverges as $\eps$ goes to $0$ because the product is singular. Wick renormalization consists in considering
$$
\mcR\PI(\xi,X)=\lim_{\eps\to0}\Big(\PI(\xi_\eps,X_\eps)-\IE\big[\PI(\xi_\eps,X_\eps)\big]\Big)
$$
in $C^{-2\kappa}(\mcS)$. This is well-known in the context of singular SPDEs, we give a proof of this convergence in Appendix \ref{SubsectinoGeometricLP}. A difference on a manifold $\CS$ compared to the torus $\IT^2$ is that the noise does not enjoy invariance in law by translation hence $\IE\big[\PI(\xi_\eps,X_\eps)\big]$ is a renormalization function and not a renormalization constant. We prove that one can consider a different choice of renormalization and actually consider some renormalization constants. In the following the parameter $z_0$ will be chosen positive, large enough depending on the size of the noise $\xi$ in our analysis. In fact it will depend on the size of the enhanced noise 
$$
\widehat{\xi} \defeq \big(\xi,\mcR\PI(\xi,X)\big)\in C^{-1-\kappa}(\CS)\times C^{-2\kappa}(\CS)
$$
thus it is important to keep track of the dependence of the renormalized product in the parameter $z_0$ as stated below. Recall that
\begin{equation}
X=(-\Delta+z_0)^{-1}(\xi)
\end{equation}
depends both on $\xi$ and $z_0\ge1$.

\medskip

\begin{proposition} \label{PropCvgce}
One has
\begin{equation*}
\sup_{\eps\in(0,1]}\Big\|\mathbb{E}\big[\PI(\xi_\eps,X_\eps)\big] - \frac{\vert \log \eps\,\vert }{4\pi}\Big\|_{C^{-2\kappa}(\CS)}<\infty
\end{equation*}
and there exists $\mcR\PI(\xi,X)\in C^{-2\kappa}(\mcS)$ such that
$$
\mcR\PI(\xi,X)=\lim_{\eps\to0}\Big({\sf \Pi}\big(\xi_\eps,X_\eps\big) - \frac{\vert \log \eps\,\vert }{4\pi}\Big)
$$
in probability in the space $C^{-2\kappa}(\mcS)$. Moreover $\mcR\PI(\xi,X)$ goes to $0$ in probability in the space $C^{-2\kappa}(\mcS)$ as $z_0$ goes to infinity.
\end{proposition}

\medskip

In order to construct the Anderson operator it remains to prove that the resolvent exists via a fixed point $R_H(z)=\Gamma(R_H(z)^\sharp)$, that is
\begin{equation}\label{EqFixedPointResolvent}
R_H(z)=-\PB_{R_H(z)}X+R_H(z)^\sharp.
\end{equation}
Recall that
\begin{align*}
&R_H(z)^\sharp =   \\
&\left(-\Delta+z_0 \right)^{-1}  \left(\text{Id}-P_\xi R_H(z)+\mcR\PI(\xi,X)R_H(z) + \overline\DC\big(R_H(z),X, \xi\big)-\PI\big(\xi,\Phi\big( R_H(z)\big)\big)+(z+z_0)R_H(z) \right) 
\\
&=\left(-\Delta+z_0 \right)^{-1}+\left(-\Delta+z_0 \right)^{-1}F(z_0,\widehat{\xi}\,)R_H(z)+(z_0+z)\left(-\Delta+z_0 \right)^{-1} R_H(z)
\end{align*}
where we used $R_H(z)^\sharp=\Phi(R_H(z))$ and 
with
\begin{equation*}
f\mapsto F(z_0,\widehat{\xi}\,)f=-\P_{\xi}f+\mcR\PI(\xi,X)f + \overline\DC(f,X,\xi)-\PI\big(\xi,\Phi(f)\big)
\end{equation*}
well-defined in the image of $R_H(z_0)$, that is for $f$ paracontrolled at first order by $X$ hence
\begin{equation}
\Phi(f)=f+\PB_fX\in H^\gamma
\end{equation}
with $\gamma>1+\kappa$ which ensures that the last term in the definition of $F(z_0,\widehat{\xi}\,)$ is well-defined. This rewrites as an equation for $R_H(z)^\sharp$, that is
\begin{equation*} \begin{split}
&R_H(z)^\sharp = \left(-\Delta+z_0 \right)^{-1}+\left(-\Delta+z_0 \right)^{-1}F(z_0,\widehat{\xi}\,) \Gamma R_H(z)^\sharp+(z_0+z)\left(-\Delta+z_0 \right)^{-1}\Gamma R_H(z)^\sharp
\end{split} \end{equation*}
hence
\begin{equation*}
\Big(\text{Id}-\left(-\Delta+z_0 \right)^{-1}F(z_0,\widehat{\xi}\,)\Gamma-(z_0+z)\left(-\Delta+z_0 \right)^{-1}\Gamma\Big)R_H(z)^\sharp= \left(-\Delta+z_0 \right)^{-1} .
\end{equation*}
The important fact is that $\left(-\Delta+z_0 \right)^{-1}F(z_0,\widehat{\xi}\,)\Gamma$ goes to $0$ in $\mcB\big(H^{1+\kappa}(\mcS), H^{2-2\kappa}(\mcS)\big)$ as $z_0$ goes to infinity. Thus we get an invertible perturbation of the identity when $z_0$ is large and the perturbation of the identity is a \textbf{compact operator} acting on $H^{1+\kappa}$. We now apply Theorem \ref{ThmFredholm} with $E=C^{-1-\kappa}(\CS)\times C^{-2\kappa}(\CS)$. Given any $\widehat{\xi_0}\in E$, there exists $z_0\ge1$ large enough such that
\begin{equation*}
R_H(z_0)^\sharp=\Big(\text{Id}-\left(-\Delta+z_0 \right)^{-1}F(z_0,\widehat{\xi}\,)\Gamma-2z_0\left(-\Delta+z_0 \right)^{-1}\Gamma\Big)^{-1}\left(-\Delta+z_0 \right)^{-1} 
\end{equation*}
exists for all $\widehat{\xi}$ such that $\|\widehat{\xi}-\widehat{\xi_0}\|_E\le1$. The theorem implies that 
\begin{equation*}
R_H(z)^\sharp=\Big(\text{Id}-\left(-\Delta+z_0 \right)^{-1}F(z_0,\widehat{\xi}\,)\Gamma-(z_0+z)\left(-\Delta+z_0 \right)^{-1}\Gamma\Big)^{-1}\left(-\Delta+z_0 \right)^{-1} 
\end{equation*}
is a well-defined meromorphic family of operators with poles of finite rank which depend continuously on $\widehat{\xi}\in C^{-1-\kappa}\times C^{-2\kappa}$. In particular, Proposition \ref{PropCvgce} states that
\begin{equation*}
\widehat{\xi_\eps}:=\Big(\xi_\eps,\PI(\xi_\eps,X_\eps)-\frac{|\log \eps|}{4\pi}\Big)
\end{equation*}
converges to $\widehat{\xi}$ in $C^{-1-\kappa}\times C^{-2\kappa}$ as $\eps$ goes to $0$ in probability. This implies the following proposition with
\begin{equation*}
R_\eps(z):=\big(-\Delta+\xi_\eps-\frac{|\log \eps|}{4\pi}-z\big)^{-1}
\end{equation*}
which is a meromorphic map of $z$ since the renormalized operator $-\Delta+\xi_\eps-\frac{|\log \eps|}{4\pi}$ is a well-defined self-adjoint operator with compact resolvent. Since $R_H(z)=\Gamma R_H^\sharp$, we get the expression
\begin{equation}
R_H(z)=\Gamma\Big(\text{Id}-\left(-\Delta+z_0 \right)^{-1}F(z_0,\widehat{\xi}\,)\Gamma-(z_0+z)\left(-\Delta+z_0 \right)^{-1}\Gamma\Big)^{-1}\left(-\Delta+z_0 \right)^{-1}.
\end{equation}


\medskip

\begin{proposition}\label{prop:imageresolvent}
The meromorphic maps $R_\eps(\cdot)$, with values in $\mcB\big(L^2(\mcS), C^{1-2\kappa}(\mcS)\big)$ converge in probability to the meromorphic map $R_H(\cdot)$ as $\eps$ goes to $0$. The map $R_H(\cdot)$ has real poles in a half-plane $\big\{\text{Re}(z)>-m\big\}$ for $m=m(\widehat{\xi})>0$ large enough.
\end{proposition}

\medskip

Now that the meromorphic map $R_H$ is constructed, one can consider the Anderson operator $H$ which corresponds to the limit
\begin{equation*}
H=\lim_{\eps\to 0}\big(-\Delta+\xi_\eps-\frac{|\log \eps|}{4\pi}\big)
\end{equation*}
which has a real non-decreasing discrete spectrum $\big(\lambda_n(\eps)\big)_{n\ge0}$.

\medskip

\begin{theorem} \label{ThmConstructionH}
The map $R_H$ is the resolvent of a closed unbounded self-adjoint operator $H$ on $L^2(\mcS)$ with real discrete spectrum bounded below. 
\end{theorem}

\medskip

\begin{proof}
Let $z_0\in\IC$ which is not a pole of the meromorphic family $R_H(\cdot)$. For $\eps_0>0$ small enough, $z_0$ is not a pole of the resolvent $R_\eps(\cdot)$ for all $\eps\in (0,\eps_0]$ hence 
$$
\lim_{\eps\to0}\|R_H(z_0)-R_\eps(z_0)\|_{\mcB(L^2\big(\CS),L^2(\CS)\big)}=0
$$ 
where $R_\eps(z_0)$ are self-adjoint operators acting on $L^2(\mcS)$. This implies that $R_H(z_0)$ itself is compact self-adjoint as an operator on $L^2(\mcS)$, consider
$$
\sigma\big(R_H(z_0)\big) = \big\{(\lambda_n-z_0)^{-1}\big\}_{n\geq 0} \subset \bbR
$$ 
with $(\lambda_n)_{n\ge0}$ increasing and $(v_n)_{n\geq 0}$ its eigenfunctions which form an orthonormal system of $L^2(\mcS)$. Since the family $R_\eps(z)$ are resolvent of operators, we get that the meromorphic family of operators $R_H(z)$ satisfies the resolvent identity
\begin{equation*}
R_H(z) = R_H(z_0)\big(\textrm{Id}+(z-z_0)R_H(z_0)\big)^{-1}
\end{equation*}
for any $z$ that is not a pole of $R_H(\cdot)$. Note that the term $\big(\textrm{Id}+(z-z_0)R_H(z_0)\big)^{-1}$ exists by meromorphic Fredholm theory in $\mcB\big(L^2(\mcS), L^2(\mcS)\big)$ relying on the compactness of $R_H(z_0)\in\mcB\big(L^2(\mcS), C^{1-2\kappa}(\mcS)\big)$. The resolvent identity implies that the range of $R_H(z_0)$ does not depend on $z_0$ thus define the $z$-independent vector space
$$
\mathfrak{D}(H) \defeq R_H(z)\big(L^2(\mcS)\big).
$$
By the resolvent equation, the meromorphic family of operators $R_H(\cdot)$ has poles contained in $(\lambda_n)_{n\geq 0}$ and satisfies the eigenvalue equation
$$ 
R_H(z)v_n = (z-\lambda_n)^{-1}v_n
$$
for $n\ge0$. This implies that we can define an unbounded operator $H$ on $L^2(\mcS)$ with domain $\mathfrak{D}(H)$ by specifying its values on an orthonormal system of $L^2(\mcS)$, that is
$$
Hv_n:=\lambda_n v_n
$$ 
for $n\ge0$. The spectrum of $H$ is bounded below since its resolvent $R_H(\cdot)$ has no poles in the half-plane $\{\text{Re}(z)\leqslant -m\}$, for $m=m(\widehat{\xi})>0$ large enough. Finally, the operator $H:\mathfrak{D}(H)\subset L^2(\mcS)\mapsto L^2(\mcS)$ is self-adjoint hence closed since
\begin{equation*}
(H-z_0)R_H(z_0) = \text{Id} : L^2(\mcS)\mapsto L^2(\mcS)
\end{equation*}
and $R_H(z_0)$ is a bounded self-adjoint operator from $L^2(\mcS)$ to $\mathfrak{D}(H)$.
\end{proof}



\section{Heat operator of the Anderson operator}
\label{SectionHeatKernelAsymptotics}

Since we have an orthonormal system $(v_n)_{n\geq 0}$ of eigenfunctions of $H$ one can consider some functions $f(H)$ of the operator $H$ for various classes of functions $f$. In this work, we are interested in the heat semigroup of $H$
$$
e^{-tH}u=\sum_{n\ge0}e^{-t\lambda_n}\langle u,v_n\rangle v_n
$$
for $u\in L^2(\CS)$. This expression defines the heat kernel as a series involving the eigenfunctions of $H$
$$
p_t(x,y)=\sum_{n\ge0}e^{-t\lambda_n}v_n(x)v_n(y),
$$
which is the main object studied in this section. In the first part we give a number of properties of the Anderson heat kernel as well as properties of $H$. In the second part we prove the main result of this work, that is a precise small time comparison of $p_t(x,y)$ with respect to the usual heat kernel $p_t^\Delta(x,y)$. In the following we refer sometimes to the eigenvalues as $\lambda_n=\lambda_n(\widehat{\xi})$ to emphasize their dependence on the enhanced noise, that is on the randomness of these eigenvalues.

\smallskip

The Anderson operator $H$ is constructed here via its resolvent
\begin{equation}
R_H(z)=\Gamma\Big(\text{Id}-\left(-\Delta+z_0 \right)^{-1}F(z_0,\widehat{\xi}\,)\Gamma-(z_0+z)\left(-\Delta+z_0 \right)^{-1}\Gamma\Big)^{-1}\left(-\Delta+z_0 \right)^{-1}
\end{equation}
hence the expression $(H-z_0)R_H(z_0)=\text{Id}$ gives
\begin{eqnarray*}
H-z_0\text{Id}&=\left(-\Delta+z_0 \right)\Big(\text{Id}-\left(-\Delta+z_0 \right)^{-1}F(z_0,\widehat{\xi}\,)\Gamma-2z_0\left(-\Delta+z_0 \right)^{-1}\Gamma\Big)\Gamma^{-1}\\
&=\left(-\Delta+z_0 \right)\Gamma^{-1}-F-2z_0
\end{eqnarray*}
thus we represent the Anderson operator as
\begin{equation}\label{ExplicitOperator}
H\Gamma u=-\Delta u+\P_{\xi}\Gamma u-\mcR\PI(\xi,X)\Gamma u - \overline\DC(\Gamma u,X,\xi)+\PI\big(\xi, u\big)
\end{equation}
for $u\in\CH^\gamma$ with $\gamma>1+\kappa$, recalling that
\begin{equation}
F(z_0,\widehat{\xi}\,)f=-\P_{\xi}f+\mcR\PI(\xi,X)f + \overline\DC(f,X,\xi) - \PI\big(\xi,\Phi(f)\big)
\end{equation}
and $\Phi=\Gamma^{-1}$. This is indeed coherent with the previous description of $H$ using the first order paracontrolled expansion from \cite{AllezChouk,GUZ,Mouzard}, and we have
\begin{equation}
\mathfrak{D}(H)\subset\Gamma(\CH^\gamma)
\end{equation}
for $\gamma>1+\kappa$. The domain $\mathfrak{D}(H)$ is explicitly described in these references where one needs strongly paracontrolled calculus or a second order paracontrolled expansion, this is not necessary in this work.


\subsection{Heat kernel and spectral properties}
\label{SubsectionHeatKernelCsq}

The solution to the parabolic Anderson model equation
\begin{equation}
\partial_tu=-\Delta u+u\xi
\end{equation}
with $u(0)=u_0$ has to be correctly renormalized as a singular SPDE. The solution is simply given by
\begin{equation*}
u(t,x)=\int_\CS p_t(x,y)u_0(y)\drm y
\end{equation*}
with the Anderson heat kernel. Even for smooth initial data, solutions are rough because of the roughness of the noise. For the heat kernel, this translates in the roughness of $p_t(x,y)$ with respect to $x,y\in\CS$ given by the following proposition.

\medskip

\begin{proposition} \label{PropKernelRegularity}
The Anderson heat kernel $p_t(x,y)$ with respect to the Riemannian volume measure on $\mcS$ is positive and continuous with respect to $(t,x,y)\in(0,\infty)\times \mcS^2$. Moreover $p_t(x,y)$ is a $(1-2\kappa)$-Hölder functions of $x,y\in\mcS$ locally uniformly in $t>0$.
\end{proposition}

\medskip

\begin{proof}
We follow the classical approach, as exposed for instance in Section 5.2 of Davies' textbook \cite{Davies}. Recall that the graph norm of $H$ on its domain $\mathfrak{D}(H)$ is defined by 
$$
\|u\|^2_H \defeq \|u\|_{L^2}^2 + \|Hu\|_{L^2}^2
$$
and that it turns $\mathfrak{D}(H)$ into a Hilbert space. With the spectral representation of the semigroup, $e^{-tH}:L^2(\CS)\to\mathfrak{D}(H)$ is continuous hence $e^{-tH}(f)\in C^{1-2\kappa}(\CS)$ for each $t>0$ and $f\in L^2(\CS)$ since we know the domain $\mathfrak{D}(H)$ is included in $C^{1-2\kappa}(\CS)$ since it is contained in the range of the resolvent $R$ and applying Proposition~\ref{prop:imageresolvent}, the reader can also refer to the proof of Theorem~\ref{ThmConstructionH}. Since $t\mapsto e^{-tH}f$ is an analytic function of $t$ on the half plane $\{\textrm{Re}(t) >0\}$ with values in the Hilbert space $\big(\mathfrak{D}(H), \|\cdot\|_H\big)$, we have that $(t,x)\mapsto e^{-tH}(f)(x)$, is a continuous function on $[t_0,t_1]\times \mcS$, for each compact interval $[t_0,t_1]\subset(0,\infty)$, analytic in the first time variable and H\"older in the second space variable. As the linear form $f\mapsto  \big(e^{-tH}f\big)(x)$, is bounded on $L^2(\mcS)$ for each $t>0$ and $x\in \mcS$ there exists $a(t,x)\in L^2(\mcS)$ such that 
$$
e^{-tH}(f)(x) = \big\langle f, a(t,x)\big\rangle_{L^2}
$$
for any $f\in L^2(\CS)$. The map 
$$
\Big( (t,x)\in(0,1]\times \mcS\Big)\mapsto a(t,x)\in L^2(\mcS),
$$ 
being weakly H\"older continuous is norm H\"older continuous with strictly smaller H\"older exponent. This is a consequence of some general principle used by Davies~\cite[Section 1.5 p.~26]{DaviesSemigroup} as follows. If we have a function $f:(t,x)\in [t_1,t_2]\times \mcS \mapsto f(t,x)\in L^2(\mcS) $ such that for all $\psi\in L^2$, $(t,x)\mapsto \left\langle f(t,x),\psi \right\rangle\in \mathcal{C}^\alpha,\alpha>0$, then $f:(t,x)\in [t_1,t_2]\times \mcS \mapsto f(t,x)\in L^2(\mcS) $ is $(\alpha-\rho)$-H\"older for all $\rho>0$. A first application of the uniform boundedness principle to the family $f(t,x)_{t\in [t_1,t_2]\times \mcS}$ which is weakly bounded in $L^2(\mcS)$ allows to deduce that $(t,x)\in [t_1,t_2]\times \mcS\mapsto f(t,x)\in L^2(\mcS)$ is strongly bounded. The family 
$$ 
\text{dist}\big((t_1,x_1),(t_2,x_2)\big)^{-\alpha} \big\vert \left\langle f(t_1,x_1)-f(t_2,x_2),\psi\right\rangle_{L^2} \big\vert 
$$
is bounded for all $\psi\in L^2(\mcS)$, then it implies by the uniform boundedness principle that
$$
\sup_{ (t_1,x_1),(t_2,x_2)} \text{dist}\big((t_1,x_1),(t_2,x_2)\big)^{-\alpha} \big( f(t_1,x_1)-f(t_2,x_2)\big) 
$$ 
is bounded in $L^2$. It follows that for all $\rho>0$ the limit as $\text{dist}((t_1,x_1),(t_2,x_2))$ goes to $0$
$$
\lim \text{dist}\big((t_1,x_1),(t_2,x_2)\big)^{-\alpha+\rho} \big( f(t_1,x_1)-f(t_2,x_2)\big) \longrightarrow 0\in L^2(\mcS),
$$ 
hence $f$ is $(\alpha-\rho)$-H\"older continuous as an $L^2(\mcS)$ valued function.

\smallskip

Then for all test functions $h_1,h_2\in C^\infty(\mcS)$, we have
\begin{align}
\big\langle e^{-tH}(h_1) , h_2\big\rangle_{L^2} &= \Big\langle e^{-\frac{t}{2}H}(h_1), e^{-\frac{t}{2}H}(h_2)\Big\rangle_{L^2} = \int p_t(x,y) h_1(x)h_2(y)\drm x\drm y
\end{align}
with 
$$
p_t(x,y) \defeq \big\langle a(t/2,x), a(t/2,y)\big\rangle_{L^2}
$$
a continuous function of its arguments. One gets the $(1-2\kappa)$-H\"older regularity of $p_t(x,y)$ as a function of $x$ for fixed $t>0$ and $y\in\CS$ since the map $(x\in\mcS)\mapsto a(t,x)\in L^2(\mcS)$ is weakly $(1-2\kappa)$-H\"older continuous it is also norm $(1-2\kappa-\rho)$-H\"older continuous for any $\rho>0$ as before. The joint regularity of $p_t(x,y)$ as a function of $(x,y)$ follows, for $0<t_0\leq t\leq t_1<\infty$.

\smallskip

Finally, the fact that $p_t(x,y)$ is positive is established in Section \ref{SubsectionConsequences} following some ideas in Cannizzaro, Friz and Gassiat in their proof of Theorem 5.1 in \cite{CFG} and our sharp description of the structure of the Schwartz kernel of $e^{-tH}$. In particular, our proof of the strong maximum principle works for all initial data in $L^2(\mcS)$ which is important to obtain a spectral gap with the argument described below.
\end{proof}

\medskip

Note that Dahlqvist, Diehl and Driver also considered in \cite{DDD} the parabolic Anderson model equation, however only with smooth initial condition hence they do not provide any insight on the heat kernel of  the Anderson operator. The next statement follows from the positivity of the heat kernel of $H$ and the Krein-Rutman theorem~\cite[Thm A.1 p.~123]{Smoller}. 

\medskip

\begin{corollary} \label{CorSimpleGroundState}
The lowest eigenvalue $\lambda_0(\widehat\xi\,)$ of $H$ is simple with a positive eigenvector almost surely.
\end{corollary}

\medskip

This question was also considered in Chouk and van Zuijlen's work \cite{CvZ}, however their proof seems incomplete since they used Cannizzaro, Friz and Gassiat's strong maximum principle \cite{CFG} which requires a continuous initial condition rather than an arbitrary initial condition in $L^2(\mcS)$. See also \cite{MouzardOuhabaz} for a simple construction of the Anderson operator which provides an elementary proof for the spectral gap based on the form domain.

\smallskip

We now state another corollary of Proposition \ref{PropKernelRegularity} that will be important for us later. It only relies on the convergence in the resolvent sense of the renormalized operators to the Anderson Hamiltonian and was already known from previous construction, see for example \cite{Labbe} if one considers only $L^2$ convergence of the ground state. For the convergence in Hölder spaces, the result is new, it could also be obtained with the description from \cite{Mouzard}. We denote as $v_{0,\eps}$ the positive ground state of the regularized operator
\begin{equation}
-\Delta+\xi_\eps-\frac{|\log \eps|}{4\pi}
\end{equation}
for $\eps>0$.

\medskip

\begin{proposition} \label{ThmContinuity}
We have the convergence in probability
\begin{equation}
\lim_{\eps\to0}\big(|\lambda_0(\widehat{\xi})-\lambda_0(\widehat{\xi_\eps})|+|\lambda_1(\widehat{\xi})-\lambda_1(\widehat{\xi_\eps})|\big)=0
\end{equation}
as well as
\begin{equation}
\lim_{\eps\to 0} \|v_0 - v_{0,\eps}\|_{C^{1-2\kappa}}=0.
\end{equation}
\end{proposition}   

\medskip

\begin{proof}
Let $\lambda\in\sigma(H)$ and $D$ be a small disc around $\lambda$ whose intersection with $\sigma(H)$ equals $\{\lambda\}$. Since the regularized and renormalized resolvent $R_\eps$ converges to $R_H$ as a Fredholm meromorphic map and $R_H(z)$ is invertible for $z\in \partial D$, we know that for $\eps$ small enough, the operators $R_\eps(z)$ are well-defined and invertible for $z\in \partial D$. Moreover it follows from the uniform convergence of $R_\eps(z)$ to $R_H(z)$ on $\partial D$ that the family of spectral projectors 
$$
\Pi_\eps^D \defeq \frac{i}{2\pi}\int_{\partial D} R_\eps(z)\drm z
$$ 
is well-defined for $\eps>0$ small enough and converges in $\mathcal{B}\big(L^2(S),H^{1-2\kappa}(\mcS)\big)$, hence the limit operator is given by
$$
\Pi^D \defeq \frac{i}{2\pi}\int_{\partial D} R_H(z)\drm z
$$
and continuous from $L^2(\mcS)$ to $H^{1-\kappa}(\mcS)$. For $z_1\notin\IR$, the meromorphic Fredholm operator $\big(\text{Id}+(z-z_1)R_\eps(z_1) \big)^{-1}$ has the same poles with multiplicity as $R_H(z)$, hence Rouché's Theorem \cite[Thm C.12]{DZ} gives that $\sigma(H_\eps)\cap D$ has fixed multiplicity for $\eps$ small enough since the poles of $R_\eps$ and $R_H$ contained in the disc $D$ have the same multiplicity. Furthermore, one has $(\Pi^D_\eps)^2 = \Pi^D_\eps$ since $\Pi_\eps^D$ is a self-adjoint spectral projector. It follows that $(\Pi^D)^2 = \Pi^D$ and that $\Pi^D$ is a self-adjoint projector such that one has
\begin{eqnarray*}
\Pi^D(v_n) = \frac{i}{2\pi}\int_{\partial D}R_H(z) (v_n) \drm z=\frac{i}{2\pi}\int_{\partial D} \big(\lambda_n-z\big)^{-1} v_n\drm z = v_n \delta_{\lambda_n=\lambda}
\end{eqnarray*}
with $v_n$ an eigenfunction of $H$. This implies that $\Pi^D$ acts as the identity when restricted on the eigenspace of $\lambda$ and vanishes on all eigenfunctions $v_n$ of eigenvalue $\lambda_n\neq \lambda$. By continuity of $\Pi^D\in \mathcal{B}\big(L^2(\mcS),L^2(\mcS)\big)$, this implies that $\Pi^D$ vanishes on the orthogonal of the eigenspace of $\lambda$ hence $\Pi^D$ is the orthogonal projector on the eigenspace of $\lambda$.

\smallskip

As a consequence of this discussion $\lambda_0(\widehat{\xi}_\eps)$ and $\lambda_1(\widehat{\xi}_\eps)$ are both converging to $\lambda_0(\widehat{\xi}\,)$ and $\lambda_1(\widehat{\xi}\,)$. By construction the lowest eigenvalues $\lambda_0(\widehat{\xi}_\eps)$ are simple for all $\eps\ge0$ however one needs a stronger result than the convergence of $\Pi_{\lambda_0(\widehat{\xi}_\eps)}$ to $\Pi_{\lambda_0(\widehat{\xi}\,)}$ in $\mathcal{B}\big(L^2(\mcS),L^2(\mcS)\big)$ to get the convergence of the ground state in $C^{1-2\kappa}(\mcS)$. Using the convergence of the kernel of $e^{-H_\eps}$ to the kernel of $e^{-H}$ in the space $\mathcal{B}\big(L^2(S),C^{1-2\kappa}(S)\big)$ which is a consequence of the continuous dependance on $\widehat{\xi}$ from Theorem \ref{ThmAsymptotic} below, we see that if one picks a small disc $D_0(\widehat{\xi}\,)$ with center $\lambda_0(\widehat{\xi}\,)$ so that $D_0(\widehat{\xi}\,)\cap \sigma(H)=\{\lambda_0(\widehat{\xi}\,)\}$, one has the convergence of 
$$
\Pi^0_\eps=e^{\lambda_0(\widehat{\xi}_\eps)} e^{-H_\eps}\Pi_\eps^{D_0(\widehat{\xi}\,)}
$$ 
to 
$$
\Pi_{\lambda_0(\widehat{\xi}\,)}= e^{\lambda_0(\widehat{\xi}\,)} e^{-H}\Pi_{\lambda_0(\widehat{\xi}\,)}
$$ 
in $\mcB\big(L^2(\mcS),C^{1-2\kappa}(\mcS)\big)$ using that $e^{H}\Pi^D=e^{\lambda}\Pi^D$. This implies the convergence of $v_{0,\eps}$ to $v_0$ in $C^{1-2\kappa}(\mcS)$. Indeed, there exists a constant $m_\eps>0$ converging to $0$ such that for all $f\in L^2(\mcS)$, one has
$$
\|\langle v_{0,\eps},f\rangle v_{0,\eps}-\langle v_0,f\rangle v_0\|_{C^{1-2\kappa}}\le m_\eps\|f\|_{L^2(\mcS)}
$$
using that the first eigenvalues are simple thus the projections are just the scalar product with the ground states. Since $(v_{0,\eps})_{\eps>0}$ is bounded in $L^2(\mcS)$, it converges weakly to $v_0'\in L^2(\mcS)$ up to an extraction. For any $z\in D_0(\widehat{\xi}\,)\backslash\{\lambda_0(\widehat{\xi}\,)\}$ and $v\in L^2(\mcS)$, we have
\begin{align*}
\big\langle (H+z)^{-1}(v_0') , f\big\rangle&=\big\langle v_0',(H+z)^{-1}(f)\big\rangle = \lim_{\eps\to0}\big\langle v_{0,\eps},(H+z)^{-1}(f)\big\rangle   \\
&= \lim_{\eps\to0}(\lambda_0(\widehat{\xi}_\eps)+z)^{-1}\big\langle v_{0,\eps},f\big\rangle = (\lambda_0(\widehat{\xi}\,)+z)^{-1}\big\langle v_0',f\big\rangle
\end{align*}
thus $v_0'=v_0$. Applying the previous bound with $f=v_0$ yields
$$
\|\langle v_{0,\eps},v_0\rangle v_{0,\eps}-v_0\|_{C^{1-2\kappa}}\le m_\eps
$$
and completes the proof. The proof shows that the spectral projectors are come continuous functions of $\widehat{\xi}$.
\end{proof}

\medskip


The following corollary states that each eigenvalue has a law absolutely continuous with respect to the Lebesgue measure. (It is not clear however that tuples of $k$ eigenvalues have a law that is absolutely continuous with respect to the Lebesgue measure in $\bbR^k$.) In particular, the first eigenvalue $\lambda_0$ has a positive probability to belong to $(-\infty,\lambda]$ for any $\lambda>0$, hence the spectrum cannot be bounded below by a deterministic constant.

\medskip

\begin{corollary} \label{CorNullKernel}
For $n\ge0$ the random variable $\lambda_n(\widehat{\xi}\,)$ has a law that is absolutely continuous with respect to the Lebesgue measure on $\bbR$, with a positive density. In particular, the kernel of $H$ is almost surely trivial and the semigroup $(e^{-tH})_{t>0}$ has no invariant Borel probability measure. 
\end{corollary}

\medskip

\begin{proof}
Given any independent random variables $X$ and $N$, the law of the random variable $X+N$ is absolutely continuous with respect to Lebesgue measure if the law of $N$ has that property. This can be seen as a regularisation of the probability density via a convolution. Thus, it suffices for example to see that the law of the random variables $\lambda_n(\widehat\xi\,)$ is absolutely continuous with respect to the law of $\lambda_n(\widehat\xi\,)+N$ with $N$ a $\mcN(0,1)$ random variable independent of $\xi$. Since the translation of the potential by a constant only induces a translation of the spectrum we have
\begin{equation*}
\lambda_n(\widehat{\xi+N}\,) = \lambda_n(\widehat\xi\,)+N
\end{equation*}
since the counterterm of the Anderson Hamiltonian for the shifted white noise $\xi+N $
does not depend on $N$ (we subtract only the singularity). It follows that the low of $\lambda_n(\widehat{\xi+N}\,)$ is absolutely continuous with respect to Lebesgue measure. The Cameron-Martin theorem gives that the law of $\xi$ is absolutely continuous with respect to the law of $\xi+N$, hence the law of $\lambda_n(\widehat{\xi}\,)$ is absolutely continuous with respect to the law of $\lambda_n(\widehat{\xi+N}\,)$ since the eigenfunctions are measurable functions of $\widehat{\xi}$ hence of $\xi$. This gives the first part of the statement.

\smallskip

Since the unbounded operator $H$ is symmetric in $L^2(\mcS)$, the heat kernel of $H$ is a symmetric function of its space arguments. So a Borel invariant probability measure has a non-negative density with respect to the Riemannian volume measure given by
$$
f(\cdot) = \int_\mcS p_t(x,\cdot)\nu(\drm x)
$$
for any $t>0$; the function is in $L^2(\mcS)$ as a continuous function on a compact set. Using the basis of eigenfunctions $(v_n)_{n\ge0}$ we get 
$$
f = \sum_{n\geq 0} c_n v_n
$$ 
and the invariance of the measure implies $e^{-tH}(f)=f$ for any $t>0$ hence $f$ is in the domain of $H$ and
$$
e^{-t\lambda_n(\widehat\xi\,)}c_n=c_n
$$
for all $n\geq 0$. The last relation implies that $f$ belongs to the kernel of $H$. Conversely, a non-null element of the kernel of $H$ defines an invariant Borel signed measure. The previous absolute continuity result implies that any eigenvalue of $H$ has null probability of being null, which finishes the proof.
\end{proof}

\subsection{Small time asymptotic for the heat kernel}
\label{SectionAsymptoticPt}



We first obtain some Schauder type estimates  for the Anderson heat semigroup via a fixed point argument and then provide a precise small time asymptotic comparison of $p_t(x,y)$ with the usual heat kernel $p_t^\Delta(x,y)$. The proofs are based on a perturbative argument hence it will be important to control
\begin{equation}
(H\Gamma+\Delta)u=\P_{\xi}\Gamma u-\mcR\PI(\xi,X)\Gamma u - \overline\DC(\Gamma u,X,\xi)+\PI\big(\xi, u\big)
\end{equation}
for $u\in H^\gamma$ with $\gamma>1+\kappa$ from expression \eqref{ExplicitOperator}. Since our goal is to study the Anderson heat semigroup it is natural to consider the conjugated operator
\begin{equation}
H^\sharp := \Gamma^{-1}H\Gamma
\end{equation}
such that
\begin{equation}
(H^\sharp+\Delta)u=\P_{\xi}\Gamma u-\mcR\PI(\xi,X)\Gamma u - \overline\DC(\Gamma u,X,\xi)+\PI\big(\xi, u\big)+\PB_{H\Gamma u}X,
\end{equation}
using that $\Gamma^{-1}=\Phi$ is explicit and we have
\begin{equation}
e^{-t H}=\Gamma^{-1}e^{-t H^\sharp}\Gamma
\end{equation}
for any $t>0$.

\medskip

\begin{lemma}\label{LemComparisonHDelta}
For any $p\in [1,\infty]$ the operator
$$
H^\sharp+\Delta: \mathcal{B}_{p,p}^{1+\kappa+\delta}(\mcS)\longmapsto \mathcal{B}_{p,p}^{-2\kappa}(\mcS)
$$
is continuous for any $\delta>0$.
\end{lemma}

\medskip

\begin{proof}
This comparison was the starting point of the construction of the Anderson operator using paracontrolled calculus \cite{GUZ,Mouzard} for $p=2$ for which the main term to control is the corrector. We have the explicit expression
\begin{equation}
(H^\sharp+\Delta)u=\P_{\xi}\Gamma u-\mcR\PI(\xi,X)\Gamma u -  \overline\DC(\Gamma u,X,\xi)+\PI\big(\xi, u\big)+\PB_{H \Gamma u}X
\end{equation}
for $u\in\CB_{p,p}^{1+\kappa+\delta}$. One needs regularity higher than $1$ to control the resonant term
\begin{equation}
\|\PI\big(\xi, u\big)\|_{\CB_{p,p}^\delta}\lesssim\|\xi\|_{\CB_{\infty,\infty}^{-1-\kappa}}\|u\|_{\CB_{p,p}^{1+\kappa+\delta}}
\end{equation}
however this is not the roughest term in the expression. Indeed the renormalized product $\mcR\PI(\xi,X)$ is of regularity $C^{-2\kappa}(\CS)$ which limits the regularity one can hope for $(H^\sharp+\Delta)u$. This is precisely the term cancelled with the strongly paracontrolled functions introduced by Allez and Chouk \cite{AllezChouk} or the second order paracontrolled expansion by the third author \cite{Mouzard} to construct the domain $\mathfrak{D}(H)$. Using the explicit expression for $\Phi=\Gamma^{-1}$, we can see that $\Gamma$ is continuous from $\mathcal{B}_{p,p}^{1+\kappa+\delta}(\mcS)$ to $\mathcal{B}_{p,p}^{1-\kappa}(\mcS)$ hence
\begin{equation}
\|\P_{\xi}\Gamma u-\mcR\PI(\xi,X)\Gamma u\|_{\CB_{p,p}^{-2\kappa}}\lesssim\|u\|_{\CB_{p,p}^{1+\kappa+\delta}}.
\end{equation}
For the corrector $\overline\DC(\Gamma u,X,\xi)$, it was proved that $v\mapsto \overline\DC(v,X,\xi)$ is continuous from $\CB_{p,p}^{1-\kappa}$ to $\CB_{p,p}^{1-3\kappa}$ for $p=2$ and $p=\infty$ thus interpolation gives the result for $p\in(2,\infty)$. For $p\in[1,2)$, one needs to adapt the proof from $p=2$. In fact, it is simple to prove that the corrector $\overline\DC(\Gamma v,X,\xi)$ belongs to $\mathcal{B}^{1-3\kappa}_{p,\infty}(\mcS)$. The last term $\PB_{H\Gamma u}X$ is controlled by all this and we get
\begin{equation}
\big\| (H^\sharp+\Delta)u \big\|_{\CB_{p,p}^{-2\kappa}}\lesssim\|u\|_{\CB_{p,p}^{1+\kappa+\delta}},
\end{equation}
which completes the proof.
\end{proof}

We can now prove the following Schauder type estimates. For the Laplace-Beltrami operator, such estimates are well-known for any regularity exponent with an explicit loss depending on the integration parameters, see for example \cite[Lem 2.6]{OhRobertTzvetkovWang}. Recall $\kappa>0$ quantifies the regularity of space white noise, by asking that $\xi\in C^{-1-\kappa}(\CS)$ almost surely.

\medskip

\begin{proposition}\label{PropSchauderH}
For any $p\in[1,\infty]$, we have
\begin{equation}
\big\| e^{-tH^\sharp}u_0 \big\|_{B_{p,p}^\alpha}\lesssim t^{-\frac{\alpha-\beta}{2}}\|u_0\|_{\CB_{p,p}^\beta}
\end{equation}
for any $\alpha\in(1+\kappa,2)$ and $\beta\in\IR$ such that $\alpha-\beta\in(0,2)$.
\end{proposition}

\medskip

\begin{proof}
Using that $H^\sharp$ is a perturbation of the Laplacian, the mild formulation for the heat equation associated to $H^\sharp$ gives
\begin{equation}
e^{-tH^\sharp}u_0=e^{t\Delta}u_0 + \int_0^te^{(t-s)\Delta}(H^\sharp+\Delta)e^{-sH^\sharp}u_0\drm s
\end{equation}
for any $t>0$. In order to prove the Schauder estimates, we prove that the map
\begin{equation}
F(u)(t):=e^{t\Delta}u_0 + \int_0^te^{(t-s)\Delta}(H^\sharp+\Delta)u(s)\drm s
\end{equation}
is a contraction on $\mathscr{E}_T=C\big([0,T],\CB_{p,p}^\beta(\CS)\big)$ equipped with the norm
\begin{equation}
\Vert u \Vert_{\mathscr{E}_T} \defeq \sup_{t\in [0,T]} \Vert u(t,.) \Vert_{\mathcal{B}^{\beta}_{p,p}}+\sup_{t\in (0,T]} t^{\frac{{\color{black} \alpha}-\beta}{2}} \Vert u(t,.) \Vert_{\mathcal{B}^{\alpha}_{p,p}}.
\end{equation}
For the high regularity part, we have
\begin{align}
\Vert F(u)(t) \Vert_{\mathcal{B}^{\alpha}_{p,p}} &\le t^{-\frac{\alpha-\beta}{2}} \Vert u_0 \Vert_{\CB_{p,p}^\beta} + C\int_0^t \vert t-s \vert^{-\frac{\alpha+2\kappa}{2}}  \Vert\left(H^\sharp+\Delta \right)u(s) \Vert_{\mathcal{B}^{-2\kappa}_{p,p}} \drm s
\\
&\le t^{-\frac{\alpha-\beta}{2}} \Vert u \Vert_{\mathscr{E}_T} + C\int_0^t \vert t-s \vert^{-\frac{\alpha+2\kappa}{2}}   \Vert u(s) \Vert_{\mathcal{B}^{\alpha}_{p,p}} \drm s
\\
& \le t^{-\frac{\alpha-\beta}{2}} \Vert u \Vert_{\mathscr{E}_T} + C\int_0^t \vert t-s \vert^{-\frac{\alpha+2\kappa}{2}} s^{-\frac{\alpha-\beta}{2}} \Vert u \Vert_{\mathscr{E}_T} \drm s\\
& \le t^{-\frac{\alpha-\beta}{2}} \Vert u \Vert_{\mathscr{E}_T} + Ct^{-\frac{\alpha+2\kappa+\alpha-\beta}{2}+1} \Vert u \Vert_{\mathscr{E}_T}
\end{align}
using Schauder estimates for the Laplace-Beltrami operator and Lemma \ref{LemComparisonHDelta}, where the condition $\alpha>1+\kappa$ is needed, with $C>0$ a constant changing from line to line. For the low regularity part, we have
\begin{align}
\Vert F(u)(t) \Vert_{\mathcal{B}^{\beta}_{p,p}}&\le \Vert u_0 \Vert_{\CB_{p,p}^\beta} + C\int_0^t \vert t-s \vert^{-\frac{\beta+2\kappa}{2}} \, \big\Vert (H^\sharp+\Delta)u(s) \big\Vert_{\mathcal{B}^{-2\kappa}_{p,p}} \drm s
\\
&\le \Vert u \Vert_{\mathscr{E}_T} + C\int_0^t \vert t-s \vert^{-\frac{\beta+2\kappa}{2}}  \Vert u(s) \Vert_{\mathcal{B}^{\alpha}_{p,p}} \drm s
\\
&\le \Vert u \Vert_{\mathscr{E}_T} + C\int_0^t \vert t-s \vert^{-\frac{\beta+2\kappa}{2}}  s^{-\frac{\alpha-\beta}{2}} \Vert u \Vert_{{\mathscr{E}_T}} \drm s \\
&\le \Vert u \Vert_{\mathscr{E}_T} + Ct^{ -\frac{\beta+2\kappa+\alpha-\beta}{2}  +1 }\Vert u \Vert_{{\mathscr{E}_T}}
\end{align}
with the same arguments and $C>0$ a constant changing from line to line.
\end{proof}

\begin{remark}
An interpolation argument allows to overcome the condition $\alpha>1+\kappa$ in the Schauder estimates and the continuity estimates also allow to get some Schauder estimates for $e^{-tH}$ -- see Proposition $4.1$ in \cite{EMR}.
\end{remark}

\medskip

The mild formulation associated to the Laplacian for the heat equation corresponding to $H^\sharp$ is
\begin{equation}
u(t)=e^{t\Delta}u_0 + \int_0^te^{(t-s)\Delta}(H^\sharp+\Delta)u(s)\drm s
\end{equation}
and this gives a comparison between the two semigroups. Indeed, the solution is also $u(t)=e^{-tH^\sharp}u_0$ hence
\begin{equation}
e^{-tH^\sharp}u_0 - e^{t\Delta}u_0 = \int_0^te^{(t-s)\Delta}(H^\sharp+\Delta)e^{-sH^\sharp}(u_0) \, \drm s
\end{equation}
for any $t>0$. This comparison was the cornerstone of the proof of Strichartz estimates by the third author and Zachhuber \cite{MouzardZachhuber} in the case of the Schrödinger group, we use it here in the case of the heat semigroup.

\medskip

\begin{theorem}\label{ThmAsymptotic}
Let $\alpha\in(-2\kappa,2-2\kappa)$ and $\beta\in(-1+2\kappa,1+2\kappa)$. For $\delta>0$, we have
\begin{equation}
\big\| e^{-tH^\sharp}u_0 - e^{t\Delta}u_0 \big\|_{\CB_{p,p}^\alpha}\lesssim t^{-\frac{\alpha-\beta}{2}+\frac{1-4\kappa}{2}}\|u_0\|_{\CB_{p,p}^{\beta}}
\end{equation}
for any $t>0$. If moreover $\alpha\le1-\kappa$, we have
\begin{equation}
\big\| e^{-tH}u_0- e^{t\Delta}u_0\big\|_{\CB_{p,p}^{\alpha\wedge (1-\kappa)}}\lesssim t^{-\frac{\alpha-\beta}{2}+\frac{1-4\kappa}{2}}\|u_0\|_{\CB_{p,p}^{\beta}}
\end{equation}
for any $t>0$.
\end{theorem}

\medskip

\begin{proof}
We have
\begin{align*}
\big\| e^{-tH^\sharp}u_0 - e^{t\Delta}u_0\big\|_{\CB_{p,p}^\alpha}&\le\int_0^t \big\| e^{(t-s)\Delta}(H^\sharp+\Delta)e^{-sH^\sharp}u_0 \big\|_{\CB_{p,p}^\alpha}\drm s\\
&\lesssim\int_0^t|t-s|^{-\frac{\alpha+2\kappa}{2}} \, \big\| (H^\sharp+\Delta)e^{-sH^\sharp}u_0\big\|_{\CB_{p,p}^{-2\kappa}}\drm s\\
&\lesssim\int_0^t|t-s|^{-\frac{\alpha+2\kappa}{2}} \big\| e^{-sH^\sharp}u_0 \big\|_{\CB_{p,p}^{1+2\kappa}}\drm s\\
&\lesssim\int_0^t|t-s|^{-\frac{\alpha+2\kappa}{2}}s^{-\frac{1+2\kappa-\beta}{2}}\|u_0\|_{\CB_{p,p}^{\beta}}\drm s\\
&\lesssim t^{-\frac{1+\alpha-\beta+4\kappa}{2}+1}\|u_0\|_{\CB_{p,p}^{\beta}}
\end{align*}
for any $\alpha\in(-2\kappa,2-2\kappa),\beta\in(-1+2\kappa,1+2\kappa)$ and $\delta>0$ using Schauder estimates for $\Delta$, Lemma \ref{LemComparisonHDelta} and Proposition \ref{PropSchauderH}. This completes the proof of the first statement since
$$
-\frac{1+\alpha-\beta+4\kappa}{2}+1=-\frac{\alpha-\beta}{2}+\frac{1-4\kappa}{2}.
$$
For the second part, we consider $\alpha<1$. We use that
\begin{align*}
e^{-tH}u_0 &= \Gamma e^{-tH^\sharp}\Gamma^{-1}u_0   \\
&=\Gamma e^{-tH^\sharp}u_0+\Gamma e^{-tH^\sharp}\PB_{u_0}X\\
&=e^{-tH^\sharp}u_0 + (\Gamma-\text{Id})e^{-tH^\sharp}u_0+\Gamma e^{-tH^\sharp}\PB_{u_0}X
\end{align*}
hence
\begin{align}
e^{-tH}u_0 - e^{t\Delta}u_0 &= \underbrace{e^{-tH^\sharp}u_0-e^{t\Delta}u_0} + \underbrace{(\Gamma-\text{Id})e^{-tH^\sharp}u_0}+\underbrace{\Gamma e^{-tH^\sharp}\PB_{u_0}X}\\
&=:R_1(t)u_0+R_2(t)u_0+R_3(t)u_0
\end{align}
for any $t>0$ where each $R_i, i\in \{1,2,3\}$ corresponds to the obvious underbraced term. The first term is controlled by the previous bound, that is
\begin{equation}
\|R_1(t)u_0\|_{\CB_{p,p}^\alpha}\lesssim t^{-\frac{\alpha-\beta}{2}+\frac{1-4\kappa}{2}}\|u_0\|_{\CB_{p,p}^{\beta}}.
\end{equation}
For the second term, we use that $\Gamma-\text{Id}$ is a regularizing operator. Indeed, $\Gamma=\Phi^{-1}$ with $\Phi(f)=f+\PB_fX$ hence
\begin{equation}
\Gamma f=\sum_{n\ge0} (-1)^n (\PB_\bullet X)^{\circ n}f.
\end{equation}
Since $X\in C^{1-\kappa}(\CS)$ and $\alpha\le1-\kappa$, this gives
\begin{align}
\|R_2(t)u_0\|_{\CB_{p,p}^{\alpha}}&\lesssim \|e^{-tH^\sharp}u_0\|_{\CB_{p,p}^{\kappa}}\\
&\lesssim t^{-\frac{(\kappa-\beta)\wedge0}{2}}\|u_0\|_{\CB_{p,p}^{\beta}}
\end{align}
using Schauder estimates for $H^\sharp$ from Proposition \ref{PropSchauderH}. For the last term, we have
\begin{align}
\|R_3(t)u_0\|_{\CB_{p,p}^{\alpha}}&\lesssim\|e^{-tH^\sharp}\PB_{u_0}X\|_{\CB_{p,p}^{1-\kappa}}\\
&\lesssim t^{\frac{\beta\wedge0}{2}}\|\PB_{u_0}X\|_{\CB_{p,p}^{1-\kappa+\beta\wedge0}}\\
&\lesssim t^{\frac{\beta\wedge0}{2}}\|u_0\|_{\CB_{p,p}^{\beta}}
\end{align}
since $\alpha\le 1-\kappa$ and again Proposition \ref{PropSchauderH}; the proof is complete.
\end{proof}

The Anderson heat kernel is related to the solution of the parabolic Anderson equation with Dirac mass as initial data. As a corollary we get the following bound for the Anderson heat kernel and the propagation of $L^2(\CS)$ initial data in Hölder spaces.

\medskip

\begin{corollary} \label{CorAsymptoticKernel}
For any $\alpha\in(-1-3\kappa,1-3\kappa)$, we have
\begin{equation}
\sup_{y\in\CS}\ \sup_{t\in(0,1]}t^{\frac{1+\alpha+\eps+3\kappa}{2}}\|p_t(\cdot,y)-p_t^\Delta(\cdot,y)\|_{C^\alpha}<\infty
\end{equation}
for any $\eps>0$ small. Moreover we have
\begin{equation}
\sup_{t\in(0,1]}t^{\frac{\alpha+3\kappa}{2}} \big\| e^{-tH}(\varphi) - e^{t\Delta}(\varphi) \big\|_{C^\alpha}\lesssim \|\varphi\|_{L^2(\CS)}.
\end{equation}
In particular the Anderson heat kernel $p_t \in C^{1-2\kappa}(\CS^2)$ depends continuously on $\widehat{\xi}$ for any positive time $t>0$.
\end{corollary}

\medskip

\begin{proof}
We have
\begin{align*}
e^{-tH}\delta_0-e^{t\Delta}\delta_0&=\Gamma e^{-tH^\sharp}\Gamma^{-1}\delta_0-e^{t\Delta}\delta_0\\ 
&=\Gamma(e^{-tH^\sharp}-e^{t\Delta})\Gamma^{-1}\delta_0+\Gamma e^{t\Delta}\Gamma^{-1}\delta_0-e^{t\Delta}\delta_0\\
&=\Gamma(e^{-tH^\sharp}-e^{t\Delta})\Gamma^{-1}\delta_0+(\Gamma-\text{Id})e^{t\Delta}\delta_0+\Gamma e^{t\Delta}\PB_{\delta_0}X
\end{align*}
hence
\begin{align*}
\|e^{-tH}\delta_0-e^{t\Delta}\delta_0\|_{C^\alpha}&\lesssim\|\Gamma(e^{-tH^\sharp}-e^{t\Delta})\Gamma^{-1}\delta_0\|_{C^\alpha}+\|(\Gamma-\text{Id})e^{t\Delta}\delta_0\|_{C^\alpha}+\|\Gamma e^{t\Delta}\PB_{\delta_0}X\|_{C^\alpha}\\
&\lesssim\|(e^{-tH^\sharp}-e^{t\Delta})\Gamma^{-1}\delta_0\|_{C^\alpha}+\|e^{t\Delta}\delta_0\|_{C^{\alpha-1+\kappa}}+\|e^{t\Delta}\PB_{\delta_0}X\|_{C^\alpha}\\ 
&\lesssim\|(e^{-tH^\sharp}-e^{t\Delta})\Gamma^{-1}\delta_0\|_{\CB_{p,p}^{\alpha+\frac{2}{p}}}+t^{-\frac{1+\alpha+2\kappa}{2}}\|\delta_0\|_{C^{-2-\kappa}}+t^{-\frac{1+\alpha+2\kappa}{2}}\|\PB_{\delta_0}X\|_{C^{-1-2\kappa}}\\ 
&\lesssim t^{-\frac{2+\alpha+\kappa}{2}+\frac{1-3\kappa}{2}}\|\Gamma^{-1}\delta_0\|_{\CB_{p,p}^{-2(1-\frac{1}{p})-\eps}}+t^{-\frac{1+\alpha+2\kappa}{2}}\|\delta_0\|_{C^{-2-\kappa}}+t^{-\frac{1+\alpha+2\kappa}{2}}\|\delta_0\|_{C^{-2-\kappa}}\|X\|_{C^{1-\kappa}}\\
&\lesssim t^{-\frac{1+\alpha+\eps+3\kappa}{2}}+t^{-\frac{1+\alpha+2\kappa}{2}}+t^{-\frac{1+\alpha+2\kappa}{2}}
\end{align*}
with $\eps>0$ using Schauder estimates and the previous bounds for $e^{-tH^\sharp}-e^{t\Delta}$. To apply the previous result, we need the conditions
\begin{align*}
\alpha+\frac{2}{p}<2-2\kappa&\quad\iff\quad p>\frac{2}{2-\alpha-2\kappa}\\ 
-2(1-\frac{1}{p})-\eps>-1+\kappa&\quad\iff\quad p<\frac{2}{1+\kappa+\eps}
\end{align*}
hence the condition
\begin{equation}
\alpha<1-3\kappa-\eps
\end{equation}
and this completes the proof of the first result. For the second bound, we have
\begin{align*}
\|e^{-tH}\varphi-e^{t\Delta}\varphi\|_{C^\alpha}&\lesssim\|\Gamma(e^{-tH^\sharp}-e^{t\Delta})\Gamma^{-1}\varphi\|_{C^\alpha}+\|(\Gamma-\text{Id})e^{t\Delta}\varphi\|_{C^\alpha}+\|\Gamma e^{t\Delta}\PB_{\varphi}X\|_{C^\alpha}\\
&\lesssim\|\Gamma(e^{-tH^\sharp}-e^{t\Delta})\Gamma^{-1}\varphi\|_{C^\alpha}+\|e^{t\Delta}\varphi\|_{C^{\alpha-1+\kappa}}+\|e^{t\Delta}\PB_{\varphi}X\|_{C^\alpha}\\
&\lesssim\|(e^{-tH^\sharp}-e^{t\Delta})\Gamma^{-1}\varphi\|_{C^\alpha}+\|e^{t\Delta}\varphi\|_{H^{\alpha+\kappa}}+\|e^{t\Delta}\PB_{\varphi}X\|_{H^{\alpha+1}}\\ 
&\lesssim\|(e^{-tH^\sharp}-e^{t\Delta})\Gamma^{-1}\varphi\|_{H^{\alpha+1}}+t^{-\frac{\alpha+\kappa}{2}}\|\varphi\|_{L^2}+t^{-\frac{\alpha+2\kappa}{2}}\|\PB_{\varphi}X\|_{H^{1-2\kappa}}\\ 
&\lesssim t^{-\frac{\alpha+1}{2}+\frac{1-3\kappa}{2}}\|\Gamma^{-1}\varphi\|_{L^2}+t^{-\frac{\alpha+\kappa}{2}}\|\varphi\|_{L^2}+t^{-\frac{\alpha+2\kappa}{2}}\|\varphi\|_{L^2}\|X\|_{C^{1-\kappa}}\\ 
&\lesssim t^{-\frac{\alpha+3\kappa}{2}}+t^{-\frac{\alpha+\kappa}{2}}+t^{-\frac{\alpha+2\kappa}{2}}
\end{align*}
using again Schauder estimates and the previous bounds for $e^{-tH^\sharp}-e^{t\Delta}$ for $\alpha+1<2-\kappa$. The continuity of the Anderson heat kernel with respect to $\widehat{\xi}$ follows from the fact that it is a solution of a fixed point problem with a map that depends continuously on the parameter $\widehat{\xi}$.
\end{proof}

\subsection{Consequences for the heat kernel}
\label{SubsectionConsequences}

We now give three consequences of our small time asymptotics. The first statement gives a property of the kernel
$$
a_t\defeq p_t-p_t^\Delta
$$
associated to the operator 
$$
A_t:=e^{-tH}-e^{t\Delta}
$$ 
for $t>0$. Theorem \ref{ThmAsymptotic} gives some pointwise estimates on the kernel $a_t$ which imply that $A_t$ is trace class in $L^2(\CS)$.

\medskip

\begin{corollary} \label{l:regremainder}
The operator $A_t$ has a well-defined Schwartz kernel $a_t(x,y)$ such that for all $\delta>0$, there exists $\rho>0$ and $T>0$ such that 
\begin{equation} \label{EqKernelEstimateA}
\sup_{y\in\mcS} \, \sup_{0<t\leq T} \, \sup_{x_1\neq x_2} \, t^{\frac{1}{2}+\delta} \, \frac{\big\vert a_t(x_1,y) - a_t(x_2,y) \big\vert}{\vert x_1-x_2 \vert^{\rho}}  +  \sup_{0<t\leq T}\sup_{x_1\neq x_2} t^{\frac{1}{2}+\delta} \vert a_t(x_1,x_2)\vert  < \infty.
\end{equation}
The operator $A_t$ is trace class in $L^2(\mcS)$ and one has   
\begin{equation*}
\text{Tr}_{L^2}(A_t) \leq \mathcal{O}\big(t^{-\frac{1}{2}-\delta}\big)
\end{equation*}
for all $t\in(0,T]$ and $\delta>0$.
\end{corollary}

\medskip

\begin{proof}
The first claim is a consequence of Theorem~\ref{ThmAsymptotic}. The key ingredient of our proof is the notion of flat trace $\text{Tr}^\flat$ which is defined for an operator $A$ with continuous kernel as
$$
\text{Tr}^\flat(A):=\int_{\mcS}A(x,x) \, \drm x.
$$
To prove the second claim, the first step is to show that for all $t>0$ the operator $e^{-tH}$ is trace class and its $L^2$-trace coincides with its flat trace. First note that 
$$
e^{-\frac{t}{2}H} = e^{\frac{t}{2}\Delta}+A_{\frac{t}{2}}
$$ 
where the operators on the right hand side have continuous Schwartz kernel by the properties of $A$ and since $t>0$ and the heat kernel is smooth at positive times. Since $p_t(x,y)\in C^0(\mcS^2)$ one has $p_t(x,y)\in L^2(\mcS^2)$ since $\mcS$ is compact with finite volume. This implies by \cite[Thm VI.23 p.~210]{RS1} that the operator $e^{-\frac{t}{2}H}$ acting on $L^2(\mcS)$ is Hilbert-Schmidt with
$$
\text{Tr}_{L^2}\left((e^{-\frac{t}{2}H})^*e^{-\frac{t}{2}H}  \right)=\int_{\mcS\times \mcS} p_{\frac{t}{2}}(x,y)\overline{p_{\frac{t}{2}}(y,x)} \, \drm x\drm y.
$$
This implies that $e^{-tH}=e^{-\frac{t}{2}H}e^{-\frac{t}{2}H} = (e^{-\frac{t}{2}H})^*e^{-\frac{t}{2}H}$ is trace class in $L^2(\CS)$ and that $\text{Tr}_{L^2}(e^{-tH})$ is well-defined to be equal to 
\begin{align*}
\int_{\mcS\times \mcS} p_{\frac{t}{2}}(x,y) \, \overline{p_{\frac{t}{2}}(y,x)}\drm x\drm y &= \int_{\mcS\times \mcS} p_{\frac{t}{2}}(x,y) \, p_{\frac{t}{2}}(y,x) \,\drm x\drm y = \int_{\mcS}p_t(x,x) \, \drm x = \textrm{Tr}^\flat\left(e^{-tH} \right) 
\end{align*}
by the Markov property of the kernel $p_t(x,y)$ with the fact that it is real-valued. The classical heat operator $e^{t\Delta}$ for $t>0$ is also trace class with $\text{Tr}_{L^2}\left(e^{t\Delta}\right)=\text{Tr}^\flat\left(e^{t\Delta} \right)$, hence the exact same properties hold true for the difference $A_t=e^{-tH}-e^{t\Delta}$. We get
\begin{eqnarray*}
\text{Tr}_{L^2}(A_t)=\text{Tr}_{L^2}\left(e^{-tH}\right)-\text{Tr}_{L^2}\left(e^{t\Delta}\right)=
\text{Tr}^\flat\left(e^{-tH}\right)-\text{Tr}^\flat\left(e^{t\Delta}\right)=\text{Tr}^\flat(A_t)
\end{eqnarray*}
hence its $L^2$-trace coincides with its flat trace. Using the first property that 
$$ 
\sup_{y\in\mcS}\Vert A_t(\cdot,y) \Vert_{\mathcal{C}^\rho}\lesssim t^{-\frac{1}{2}-\delta },
$$ 
we conclude that 
$$
\text{Tr}_{L^2}(A_t)=\int_{\mcS} A_t(x,x)\drm x=\mathcal{O}\big( t^{-\frac{1}{2}-\delta }\big)
$$ 
which is the desired claim.
\end{proof}

We prove that the strong maximum principle for the semigroup $e^{-tH}$ follows from our method of proof of Theorem \ref{ThmAsymptotic}. We follow Cannizzaro, Friz and Gassiat's proof \cite{CFG}.

\medskip

\begin{proposition}\label{prop:max}
For any non-negative $u_0\in L^2$ and $t>0$, the function $e^{-tH}(u_0)$ is continuous and has a positive minimum.
\end{proposition}

\medskip

\begin{proof}
For $u_0\in L^2(\CS)$, we have $e^{-tH}(u_0) \in\mathcal{D}(H)$ hence it belongs to $C^{1-\kappa}(\CS)$. The Feynamn-Kac formula for the operator $H_\eps=-\Delta+\xi_\eps-c_\eps$ gives
\begin{equation}
e^{-tH_\eps}(u_0)(x) = \mathbb{E}\Big[ e^{-2tc_\eps -\int_0^{t}2\xi_\eps(B_s)\drm s }u_0(B_t^x) \Big] > 0  
\end{equation}
with $(B_t^x)_{t\ge0}$ a Brownian motion starting at $x\in\CS$ however the limit when $\eps$ goes to $0$ only yields $e^{-tH}(u_0)\ge0$. In order to prove that $e^{-tH}u_0>0$, we prove that the Anderson heat kernel is positive $p_t(x,y)>0$ for all $x,y\in\CS$.

\smallskip

Let $D$ be the diameter of the surface $\mcS$. On the one hand, Corollary \ref{CorAsymptoticKernel} gives
\begin{equation}
p_t(x,y)=p_t^\Delta(x,y)+t^{-\frac{1+\delta}{2}}R_H(x,y)
\end{equation}
with $r\in L^\infty(\CS^2)$. On the other hand, by the Li-Yau estimates \cite[Thm 4.8 p.~172]{SchoenYau}, we have a lower bound of the form
$$
p_t^\Delta(x,y) \ge\frac{c_1}{t} e^{-\frac{c_2d^2(x,y)}{t}}
$$ 
with $c_1,c_2>0$ positive constants. We get
\begin{align}
p_t(x,y)&\ge \frac{c_1}{t}e^{-\frac{c_2D^2}{t}}-\frac{c_3}{t^{\frac{1+\delta}{2}}}\\
&\ge\frac{1}{t}(c_1e^{-\frac{c_2D^2}{t}}-c_3t^{\frac{1-\delta}{2}})
\end{align}
with $c_3=\sup_{x,y\in\CS}|R_H(x,y)|>0$. For $t\le t_0$ with $t_0$ small enough depending on $D,c_1,c_2,c_3>0$, the right hand side is positive. This completes the proof since
\begin{equation}
p_t(x,y)=\int_\CS p_{t-t_0}(x,z)p_{t_0}(z,y)\drm z
\end{equation}
for $t>t_0$ and $p_{t-t_0}(x,z)\ge0$
\end{proof}

\smallskip

Recall we denote by $(v_n)_{n\geq 0}$ an orthonormal system of $L^2$ normalized eigenfunctions of $H$, with $Hv_n=\lambda_n v_n$ and the $\lambda_n$ ordered in non-decreasing order. Finally one can consider the ground state transform of the Anderson operator, that is
\begin{equation}
\SL=v_0^{-1}(H-\lambda_0)v_0.
\end{equation}
Since $v_0$ is a continuous positive function on the compact manifold $\CS$, we have $v_0^{-1}\in L^\infty(\CS)$. Indeed, the ground state $v_0$ is non-negative as the $L^\infty$ limit of the positive functions $v_{0,\varepsilon}$ and it satisfies 
$$
e^{t\lambda_0}e^{-tH}(v_0) = v_0 
$$ 
for all $t> 0$ and the strong maximum principle implies the positivity of the ground state. It follows that the Anderson heat semigroup satisfies Gaussian upper and lower bound as well as similar bounds for its Green function, see Stroock's book \cite{StroockPDEs} for a reference in the context of a smooth potential. For any $c\ge-\lambda_0(\widehat{\xi})$, the resolvent $R_H(c)=(H-c)^{-1}$ is well defined with kernel $G_c$, we consider the Green function $G=G_c$ for fixed such random constant in the following. We also denote as $G^\Delta$ the Green function of the Laplace-Beltrami operator on $\CS$.

\medskip

\begin{proposition} \label{PropHeatKernalBounds}
For any $T>0$, there exists random constants $c,m>0$ such that
\begin{equation}
m^{-1}p_{c^{-1}t}^\Delta(x,y)\le p_t(x,y)\le mp_{ct}^\Delta(x,y)
\end{equation}
for any $t\in(0,T]$ and $x,y\in\CS$. Moreoever, there exists a random constant $\ell>0$ such that
\begin{equation}
\ell^{-1}G^\Delta(x,y)\le G(x,y)\le \ell G^\Delta(x,y)
\end{equation}
for any $x,y\in\CS$. In particular, the heat semigroup $(e^{-tH})_{t>0}$ is hypercontractive. 
\end{proposition}

\medskip

\begin{proof}
We prove two-sided Gaussian bounds for the heat kernel $p_t^\eps(x,y)$ associated to 
\begin{equation}
H_\eps=-\Delta+\xi_\eps-c_\eps
\end{equation}
uniform with respect to $\eps>0$ for fixed $t>0$. Since $p_t$ depends continuously on $\widehat{\xi}$ in $C^{1-\kappa}(\CS^2)$ with Corollary \ref{CorAsymptoticKernel}, the result will follow from the convergence of $\widehat{\xi_\eps}$ to $\widehat{\xi}$. We consider the ground state transform of the operator, that is
\begin{equation}
\SL_\eps:=v_{0,\eps}^{-1}(H_\eps-\lambda_{0,\eps})v_{0,\eps}
\end{equation}
with $v_{0,\eps}>0$ the ground state associated to $H_\eps$ and $\lambda_{0,\eps}=\lambda_0(\widehat{\xi_\eps})$. This is a conservative perturbation of the Laplace-Beltrami operator for which we have two-sided Gaussian bounds on the heat kernel $p_t^{\SL_\eps}$ using that $v_{0,\eps}$ is continuous, see for instance Section 4.3 and Section 6.4 of Stroock's book \cite{StroockPDEs}.
So
there is a continuous positive function $c(\cdot)$
with $c(0)=1$ such that setting 
$c_\eps=c\big(\text{osc}(\log v_{0,\eps})\big)$ where the oscillations of $\log v_{0,\eps}$ are defined as $\text{osc}(\log v_{0,\eps}):=\sup_\CS \log v_{0,\eps}-\inf_\CS \log v_{0,\eps}$ for any function $\log v_{0,\eps}$,
one has
\begin{equation}
\frac{1}{c_\eps t}\,\exp\left(-\frac{c_\eps d(x,y)^2}{t}\right) \leq p_t^{\SL_\eps}(x,y) \leq \frac{c_\eps}{t}\,\exp\left(-\frac{d(x,y)^2}{c_\eps t}\right)
\end{equation}
for all $0<t\leq 1$ and $x,y\in\mcS$. Since
\begin{equation}
p_t^\eps(x,y)=e^{t\lambda_{0,\eps}}\frac{v_{0,\eps}(x)}{v_{0,\eps}(y)}p_t^{\SL_\eps}(x,y)
\end{equation}
for any $t>0$, we get
\begin{equation}
\frac{e^{-t\lambda_{0,\eps}}}{m_\eps c_\eps t}\,\exp\left(-\frac{c_\eps d(x,y)^2}{t}\right) \leq p_t^\eps(x,y) \leq \frac{m_\eps c_\eps e^{-t\lambda_{0,\eps}}}{t}\,\exp\left(-\frac{d(x,y)^2}{c_\eps t}\right)
\end{equation}
with $m_\eps=\frac{\max_\CS v_{0,\eps}}{\min_\CS v_{0,\eps}}>0$. To conclude, one only needs to prove that the families $(c_\eps)_{\eps>0}$ and $(m_\eps)_{\eps>0}$ are bounded from below by a constants $c,m>0$ since $\lambda_{0,\eps}$ converges to $\lambda_0\in\IR$ with Proposition \ref{ThmContinuity}. The same proposition gives that $v_{0,\eps}$ converges to $v_0$ in $C^{1-2\kappa}$. While Krein-Rutman Theorem only gives $v_0>0$ almost everywhere on $\CS$ for the measure $\mu$, the strong maximum principle from Proposition \ref{prop:max} gives that $v_0$ has a minimum strictly positive on $\CS$ which completes the proof of the two-sided Gaussian bounds. The estimates on the Anderson Green functions follow directly from the two-sided Gaussian bounds and the expression
\begin{equation}
(H+c)^{-1}=\int_0^\infty e^{-t(H+c)}\drm t
\end{equation}
as well as the same expression for the Laplace-Beltrami operator.
\end{proof}

\subsection{Consequences for the spectrum}
\label{SubsectionEigenBounds}

In this section we prove a number of results on the eigenvalues and eigenfunctions of $H$ using the sharp description of $p_t$ given by Theorem \ref{ThmAsymptotic}. The first statement is that
\begin{eqnarray} \label{EqTraceEstimate}
\text{Tr}_{L^2}(e^{-tH}) = \text{Tr}_{L^2}(e^{t\Delta}) + \text{Tr}_{L^2}(A_t) = \frac{\mu(\mcS)}{4\pi t} + \mathcal{O}\big(t^{-\frac{1}{2}-\delta}\big).
\end{eqnarray}
for any $\delta>0$ which allows to recover Weyl's law for the Anderson operator initially proved via a minimax representation of the eigenvalues in \cite{Mouzard}. If folllows here from the small time equivalent for the heat kernel by Karamata's Tauberian Theorem~\cite[Thm 2.42 p.~94]{BGV}.

\medskip

\begin{proposition} \label{CorWeyl}
We have almost surely the equivalent
\begin{eqnarray}
\Big|\big\{\lambda\in\sigma(H)\,;\,\lambda\leq a\big\}\Big|\underset{a\to+\infty}{\sim} \frac{\mu(\mcS)}{4\pi}\,a.
\end{eqnarray}
\end{proposition}

\medskip

One thus has almost surely the equivalent
$$
\lambda_n(\widehat\xi\,)\underset{n\to\infty}{\sim}\frac{4\pi}{\mu(\mcS)}\,n
$$ 
as $n$ goes to $\infty$, the same asymptotics as the eigenvalues of the Laplace-Beltrami operator. While previous works provided some bounds on the tails of the eigenvalues, we provide here some upper bound on the density. Description of the tails of the eigenvalues were given in some previous works such as \cite{AllezChouk,Labbe,Mouzard}.

\medskip

\begin{proposition}
One has
$$
\bbP\big(1\leq \vert\lambda_k(\widehat{\xi}\,)\vert \leq \lambda\big) \lesssim \Big(\frac{\lambda}{k}\Big)^b
$$
for all $k\geq 1$ and $\lambda\geq 1$.
\end{proposition}

\medskip

\begin{proof}
Since $\lambda_n\ge\lambda_k$ for $k\le n$, we get
\begin{equation}
ne^{-t\lambda_n}\le \sum_{n\ge0}e^{-t\lambda_n}
\end{equation}
hence taking $t=|\lambda_n^{-1}|$ in the bound
\begin{equation}
\text{Tr}_{L^2}(e^{-tH}) \leq \frac{c_1(\widehat{\xi}\,)}{t}
\end{equation}
for $t\in(0,1]$ yields
\begin{equation}
\big\vert\lambda_k(\widehat{\xi}\,)\big\vert \geq \frac{e}{c_1(\widehat{\xi}\,)}\,k,
\end{equation}
conditionned to the fact that $|\lambda_n|\ge1$ to ensure $t\in(0,1]$. The function 
$$
F_1(x) \defeq \bbP\big(c_1(\widehat{\xi}\,)\geq x\big)
$$ 
has thus the property that 
$$
\bbP\big(1\leq \vert\lambda_k(\widehat{\xi}\,)\vert \leq \lambda\big) \leq F_1\Big(\frac{ek}{\lambda}\Big)
$$
for all $k\geq 1$ and $a\geq 1$. The analysis of the proof of Theorem \ref{ThmAsymptotic} shows that one can choose  $c_1(\widehat{\xi}\,)$ of the form
$$
c_1(\widehat{\xi}\,) = e^{c\Vert\widehat{\xi}\,\Vert_{C^{-1-\kappa}(\CS)\times C^{1-2\kappa}(\CS)}},
$$
for a positive constant $c$. As we know that $\xi$ has a Gaussian tail and $\CR{\sf \Pi}(\xi,X)$ has an exponential tail, see e.g. Proposition 2.2 in \cite{Mouzard}, there exists a positive constant $b$ such that 
$$
F_1(x)\lesssim \frac{1}{x^b}.
$$
\end{proof}

We also get some bounds on the growth of the Hölder norms of the eigenfunctions of the Anderson operator. Similar bounds were obtained in the $L^q$ spaces in \cite{MouzardZachhuber} which are sharper than the one we get here. However, we are able to consider Hölder spaces since we work with the heat operator instead of the Schrödinger semigroup. 

\medskip

\begin{proposition} \label{ThmEstimateEigenfunctions}
For any $\alpha\in(0,1)$ there exists a random constant $C>0$ such that
\begin{equation} \label{EqHolderNormEn}
\Vert v_n\Vert_{C^\alpha} \le C\big(1+\lambda_n^{\frac{\alpha+1+\delta}{2}}\big)
\end{equation}
for any $\delta>0$.
\end{proposition}

\medskip

\begin{proof}
We have
\begin{equation}
e^{-t\lambda_n}v_n=e^{t\Delta}v_n+A_tv_n
\end{equation}
with $A_t=e^{-tH}-e^{t\Delta}$. For $\alpha\in(-1-2\kappa,1-2\kappa)$, Corollary \ref{CorAsymptoticKernel} gives
\begin{equation}
\big\| A_t(v_n) \big\|_{C^\alpha}\lesssim t^{-\frac{\alpha+\delta+3\kappa}{2}}
\end{equation}
for any $\delta>0$ since $\|v_n\|_{L^2}=1$ and Schauder estimates for $\Delta$ gives
\begin{align}
\big\| e^{t\Delta}(v_n) \big\|_{C^\alpha}&\lesssim t^{-\frac{\alpha+1+\delta}{2}}\|v_n\|_{C^{-1-\delta}} \lesssim t^{-\frac{\alpha+1+\delta}{2}}\|v_n\|_{H^{-\delta}} \lesssim t^{-\frac{\alpha+1+\delta}{2}}
\end{align}
for any $\delta>0$. We get
\begin{equation}
e^{-t\lambda_n}\|v_n\|_{C^\alpha}\lesssim t^{-\frac{\alpha+1+\delta}{2}}+t^{-\frac{\alpha+\delta+3\kappa}{2}}
\end{equation}
hence taking $t=\lambda_n^{-1}$ with belongs to $(0,1]$ for $n$ large enough gives
\begin{equation}
\|v_n\|_{C^\alpha}\lesssim \lambda_n^{\frac{\alpha+1+\delta}{2}}+\lambda_n^{\frac{\alpha+\delta+3\kappa}{2}}
\end{equation}
which completes the proof. Note that this also follows directly for the general Schauder estimates given by Proposition $4.1$ in \cite{EMR}.
\end{proof}

We conclude this section with two lower bounds on the spectral gap of $H$ under two kinds of assumptions, geometric and functional analytic. For any smooth volume measure $\nu$ on $\mcS$ one can define its Cheeger constant
$$
C(\nu) \defeq \underset{A\subset\mcS}{\inf} \ \frac{\sigma_\nu(\partial A)}{\min\big\{\nu(A),\nu(\mcS\backslash A)\big\}}
$$
where
$$
\sigma_\nu(\partial A) \defeq \underset{\kappa\searrow 0}{\liminf} \ \frac{\nu(A^{(\kappa)}) - \nu(A)}{\kappa}
$$
with $A^{(\kappa)} \defeq \{m\in\mcS\,;\,d(m,A)\leq \kappa\}$ the $\kappa$-enlargement of a set $A\subset\CS$. Recall we denote by $\mu$ the Riemannian volume measure on $\mcS$.

\medskip

\begin{proposition} \label{ThmSpectralGapCheeger}
One has almost surely the following estimate on the spectral gap
$$
\lambda_1(\widehat{\xi}\,) - \lambda_0(\widehat{\xi}\,) \geq \frac{C(v_0^2\mu)^2}{4}
$$
with $u_0$ the Anderson ground state.
\end{proposition}

\medskip

\begin{proof}
Proceeding as in the proof of Proposition \ref{PropHeatKernalBounds}, we see that it suffices to prove that the spectral gap $\lambda_1(\widehat{\xi}_\eps) - \lambda_0(\widehat{\xi}_\eps)$ of the conjugated regularized operator 
$$
-\Delta - 2(\nabla  \log v_{0,\eps})\nabla
$$ 
is bounded below by $C(v_{0,\eps}^2\mu)^2/4$. Indeed, the convergence of $v_{0,\eps}$ to $u_0$ in $C^{1-\kappa}(\mcS)$ proved in Proposition \ref{ThmContinuity} implies that $C(v_{0,\eps}^2\mu)$ is converging to $C(v_0^2\mu)$ as $\eps$ goes to $0$. While it is classical in Riemannian geometry, we prove the Cheeger lower bound on $\lambda_1(\widehat{\xi}_\eps) - \lambda_0(\widehat{\xi}_\eps)$ adapted to our context. We consider the measure 
$$
\nu_{0,\eps} := v_{0,\eps}^2\mu.
$$
For all smooth functions $f\in C^\infty(\mcS)$ with median value $m_{0,\eps}(f)$ with respect to $\nu_{0,\eps}$, one has
\begin{equation} \label{EqGradientEstimate}
\int_S \Vert \nabla f \Vert\, d\nu_{0,\eps} \geqslant C(\nu_{0,\eps})\int_{\mcS} \big\vert f-m_{0,\eps}(f)\big\vert \, d\nu_{0,\eps}.
\end{equation} 
If one takes \eqref{EqGradientEstimate} for granted for a moment, one can apply this inequality to the function $f\vert f\vert$ where $f$ is rescaled in such a way that it has unit $L^2(\nu_{0,\eps})$-norm and $f^{-1}(0)$ and $\left(f\vert f\vert\right)^{-1}(0) $ have equal $\nu_{0,\eps}$-measure $\nu_{0,\eps}(\mcS)/2$, so $f\vert f\vert$ has a null median. This yields
\begin{eqnarray*}
\int_{\mcS} \big\Vert \nabla \left( f\vert f\vert\right)\big\Vert \, d\nu_{0,\eps} = 2\int_{\mcS} \Vert f \nabla f\Vert\, d\nu_{0,\eps} \geqslant C(\nu_{0,\eps})\int_{\mcS} \vert f\vert ^2\,d\nu_{0,\eps} = C(\nu_{0,\eps})
\end{eqnarray*}
and we get from Cauchy-Schwartz inequality that
\begin{eqnarray*}
C(\nu_{0,\eps}) \leq  2\Vert \nabla f \Vert_{L^2(\nu_{0,\eps})}.
\end{eqnarray*}
In the general case if $f\in C^\infty(S,\mathbb{R})$ is such that $\int_\mcS fd\nu_{0,\eps}=0$ and $\int_S f^2d\nu_{0,\eps}=1$, one can use the inequality 
$$
\int_\mcS (f+c)^2d\nu_{0,\eps} =\int_\mcS (f^2+c^2)\,d\nu_{0,\eps}\geqslant \int_S f^2\,d\nu_{0,\eps}
$$ 
to possibly add a constant to $f$ and trade the assumption that $\int_\mcS fd\nu_{0,\eps}=0$ for the assumption that $f^{-1}(0)$ cuts $\mcS$ in two pieces of equal $\nu_{0,\eps}$ measure. Applying the above arguments to $\frac{f+c}{\Vert f+c \Vert_{L^2(\nu_{0,\eps})}}$ yields
\begin{eqnarray*}
C(\nu_{0,\eps}) \leq  2 \frac{ \Vert \nabla f \Vert_{L^2(\nu_{0,\eps})} }{\Vert f+c \Vert_{L^2(\nu_{0,\eps})}}\leq 2 \frac{ \Vert \nabla f \Vert_{L^2(\nu_{0,\eps})} }{\Vert f\Vert_{L^2(\nu_{0,\eps})}}.
\end{eqnarray*}
The representation of the spectral gap of $\Delta+2(\nabla \log v_{0,\eps})\nabla$ as a Rayleigh quotient
$$
\lambda_1(\widehat{\xi}_\eps) - \lambda_0(\widehat{\xi}_\eps) = \inf_{\int_{\mcS} f\,d\nu_{0,\eps} = 0}\, \frac{\int_{\mcS} \Vert \nabla f\Vert^2\,d\nu_{0,\eps}}{\int_{\mcS} \vert f\vert^2\,d\nu_{0,\eps}}
$$
gives
$$
\lambda_1(\widehat{\xi}_\eps) - \lambda_0(\widehat{\xi}_\eps) \geq \frac{C(\nu_{0,\eps})^2}{4}.
$$
It remains to prove formula \eqref{EqGradientEstimate}. Recall from the coarea formula that one has 
$$
\int_{\mcS} \Vert \nabla f\Vert\,d\nu_{0,\eps}  =\int_{\bbR} \sigma_{\nu_{0,\eps}}\big(\{f=t\}\big)\,dt.
$$ 
From the isoperimetric inequality
$$
 \sigma_{\nu_{0,\eps}}(\partial A) \geqslant C(\nu_{0,\eps}) \min\big(\nu_{0,\eps}(A),\nu_{0,\eps}(\mcS\setminus A)\big)
$$ 
we deduce that if $0$ is a median of $f$ we have the bounds
\begin{align}
\int_{\mcS} \Vert \nabla f\Vert\,d\nu_{0,\eps} &= \int_{f\leq 0} \vert \nabla f\vert\,d\nu_{0,\eps} + \int_{f>0} \vert \nabla f\vert\,d\nu_{0,\eps}  = \int_{-\infty}^0 \sigma_{\nu_{0,r}}\left(\{f=t\}\right)  dt+\int_{0}^\infty \sigma_{\nu_{0,\eps}}\left(\{f=t\}\right) \, dt   \\
&\geqslant C(\nu_{0,\eps}) \left(\int_{-\infty}^0 \nu_{0,\eps}(\{f\leq t\}) dt + \int_{0}^\infty \nu_{0,\eps}(\{f>t\}) \, dt \right) \geqslant C(\nu_{0,\eps}) \int_{\mcS} \vert f\vert\,d\nu_{0,\eps}
\end{align}
where we used integration by parts for the last step and disintegration of the volume $\nu_{0,\eps}$ along level sets of $f$. 
\end{proof}   

For a non-negative measure $\nu$ on $\mcS$, one can consider the $\nu$-entropy of a positive integrable function $f$ such that $\int_\mcS f\vert \log f\vert\,\drm\nu<\infty$ as the  quantity
$$
\textrm{Ent}_\nu(f) \defeq \int_\mcS f \log f\,\drm\nu - \left(\int_\mcS f\,\drm\nu\right) \log\left(\int_\mcS f\,\drm\nu\right).
$$
We say that a measure $\nu$ on $\mcS$ satisfies a $\log$-Sobolev inequality with constant $C_{{\textsc{LS}}}$ with respect to the Dirichlet form associated with the Riemannian gradient operator $\nabla$ if
$$
\textrm{Ent}_\nu(f^2) \leq 2C_{{\textsc{LS}}} \int_\mcS\vert\nabla f\vert^2\,\drm\nu
$$
for all functions $f$ in the domain of the Dirichlet form. Such an inequality is known to imply a Poincaré inequality with constant $C_{{\textsc{LS}}}^{-1}$ and a corresponding spectral gap. Bakry, Gentil and Ledoux's monograph \cite{BakryGentilLedoux} presents several geometric conditions ensuring that $\mu$ satisfies a $\log$-Sobolev inequality.   

\medskip

\begin{proposition} \label{ThmSpectralGapH}
Assume that the Riemannian volume form $\mu$ satisfies a $\log$-Sobolev inequality with constant $C_{{ \text{LS}}}$. Then the spectral gap of $H$ satisfies almost surely the lower bound
$$
\lambda_1(\widehat{\xi}\,) - \lambda_0(\widehat{\xi}\,) \geq \bigg(\frac{\min v_0}{\max v_0}\bigg)^2 \, \frac{\big(\max v_0^4+\max v_0^{-4}\big)^{-1}}{2C_{{ \text{LS}}}}
$$
with $u_0$ the Anderson ground state. 
\end{proposition}

\medskip

\begin{proof}
Denote by $m_\eps$ the spectral gap of $H_\eps$ in $L^2(\mu)$ and by $m_\eps'$ the spectral gap of $H_\eps$ in $L^2(v_{0,\eps}^{-2}\mu)$. Then $m_\eps'$ is equal to the spectral gap of the conjugated operator $\Delta - 2\nabla(\log v_{0,\eps})\nabla$ and
$$
m_\eps \geq m_\eps' \, \bigg(\frac{\min v_{0,\eps}}{\max v_{0,\eps}}\bigg)^2.
$$
As in the proof of Theorem \ref{ThmSpectralGapCheeger}, we recognize in the conjugated operator the Dirichlet form of the Riemannian gradient operator with respect to the weighted Riemannian volume form $v_{0,\eps}^2\mu$. As Holley and Stroock stability argument for $\log$-Sobolev inequality ensures that the weighted measure $v_{0,\eps}^2\mu$ satisfies, under the assumption of the statement, a $\log$-Sobolev inequality with constant $2C_{{ \textsc{LS}}} \, \big(\max v_{0,\eps}^4+\max v_{0,\eps}^{-4}\big)$, we see that
$$
m_\eps' \geq \frac{\big(\max v_{0,\eps}^4+\max v_{0,\eps}^{-4}\big)^{-1}}{2C_{{ \textsc{LS}}}},
$$
see for example Proposition 5.1.6 in \cite{BakryGentilLedoux} for a proof of the stability argument. We thus have the lower bound
$$
\lambda_{1,\eps} - \lambda_{0,\eps} = m_\eps \geq \bigg(\frac{\min v_{0,\eps}}{\max v_{0,\eps}}\bigg)^2 \, \frac{\big(\max v_{0,\eps}^4+\max v_{0,\eps}^{-4}\big)^{-1}}{2C_{{ \textsc{LS}}}}
$$
and we conclude by using the continuity of the eigenvalues as functions of $\widehat{\xi}_\eps$ and the convergence in $L^\infty(\mcS)$ of $v_{0,\eps}$ to $v_0$ from Proposition \ref{ThmContinuity}.
\end{proof}

\section{Anderson Gaussian free field}
\label{SectionAndersonGFF}


We fix throughout this section a random variable 
$$ 
c > -\lambda_0(\widehat\xi\,)
$$
such that the operator $H+c$ is positive and defines a distribution-valued Gaussian field with covariance $(H+c)^{-1}$. We call it the Anderson Gaussian free field $\phi$ and it can be defined by the formula
\begin{equation}
\phi \defeq \sum_{n\geq 0} \frac{\gamma_n}{(c+\lambda_n)^{\frac{1}{2}}} \, v_n
\end{equation}
where the $\gamma_n$ are independent, identically distributed, real-valued random variables with law $\mcN(0,1)$. Note that this is a random field $\phi$ in a random environment $\xi$ hence it has two independent layers of randomness, one coming from $H$ and the other coming from the $\gamma_n$, a notation emphasizing that fact would be
$$
\phi(\omega,\omega') = \sum_{n\geq 0} \frac{\gamma_n(\omega')}{(c(\omega)+\lambda_n(\omega))^{\frac{1}{2}}} \, v_n(\omega).
$$
In the following, we will only take expectation with respect to this new environement conditionned on the random environment $\xi$. With this in mind, the random field $\phi$ is a centered Gaussian field with covariance is given by
\begin{equation}
\IE\big[\phi(x)\phi(y)\big]=G(x,y)
\end{equation}
for $x,y\in\CS$, we do not keep the dependence with respect to $c$ to lighten the notation which ensure that $G$ is positive. We refer to Da Prato's book \cite{DP} for general results on Gaussian measure in Hilbert spaces. We proved in Proposition \ref{PropHeatKernalBounds} that the Green function has a logarithmic divergence along the diagonal thus the Anderson Gaussian Free Field is a log correlated Gaussian fields and takes values in distribution. The Kolmogorov criterion immediatly gives the following result.

\medskip

\begin{proposition}
The Anderson Gaussian free field is almost surely in $C^{-\delta}(\mcS)$ for every $\delta>0$.
\end{proposition}

\medskip

A natural space to consider is the associated Cameron-Martin space
\begin{equation}
\CC\CM=(H+c)^{-\frac{1}{2}}L^2(\CS)
\end{equation}
associated to the Anderson Gaussian free field, see Section 1.7 in \cite{DP} for details. One motivation behind this space is that the law of $\phi+f$ is absolutely continuous with respect to the law of $\phi$ if and only if $f\in\CC\CM$, while this is always the case in finite dimension. In particular we have
\begin{equation}
\CC\CM\subset H^{1-\kappa}(\CS)
\end{equation}
from our construction of the Anderson operator. This space will appear in the proof of the following result. For $n\geq 2$ consider
$$
a_n:=\int_{\CS^n} \prod_{i=1}^n G(x_i,x_{i+1}) \drm x_1\dots\drm x_n
$$
with the convention $x_{n+1}=x_1$ in the integral. Since the Green function has a logarithmic divergence near the diagonal with Proposition \ref{PropHeatKernalBounds}, this is indeed well-defined. We have
$$
a_n = \text{Tr}_{L^2}\big((H+c)^{-n}\big)
$$ 
hence the quantity $a_n$ is purely spectral as we have
\begin{equation} \label{EqTraceFormulaAn}
a_n = \sum_{k\geq 0} \big(\lambda_k(\widehat\xi\,)+c\big)^{-n}
\end{equation}
from Lidskii's theorem. Consider the heat regularized Anderson Gaussian free field
$$
\phi_\eps = e^{\eps\Delta}(\phi)
$$
for $\eps>0$. We define its regularized Wick square as
$$
\Wick{\phi_\eps^2}\ := \phi_\eps^2 - \IE\big[\phi_\eps^2\big]
$$
where $\IE\big[\phi_\eps^2\big]$ is a divergent quantity. While $(H+c)^{-1}$ is not trace class, it will be crucial in the proof of the next statement that $(H+c)^{-1}$ is Hilbert-Schmidt, which is ensured by the Weyl law from Corollary \ref{CorWeyl}. We consider the partition function
\begin{equation}
Z(\lambda):=\IE\Big[e^{-\lambda\langle\Wick{\phi^2},1\rangle}\Big]
\end{equation}
for $\lambda\in\IC$.

\medskip

\begin{theorem} \label{ThmWickCharacterizationSpectrumH}
There exists a random distribution $\Wick{\phi^2}$ such that for any $\delta>0$, we have
\begin{equation}
\lim_{\eps\to0}\|\Wick{\phi^2}-\Wick{\phi_\eps^2}\|_{C^{-\delta}}=0
\end{equation}
in probability. For all $\lambda\in \bbC$ sufficiently small one has
\begin{equation} \label{EqPartitionFunction}
Z(\lambda) = {\det}_2\Big(\textrm{Id} + \lambda (H+c)^{-1}\Big)^{-1/2} = \exp\left(\sum_{n\geq 2}\frac{(-\lambda)^n a_n}{2n}\right).
\end{equation}
Moreover this function of $\lambda$ has an analytic extension to all of $\bbC$.
\end{theorem}

\medskip

\begin{proof}
Proposition 9.3.1 in Glimm and Jaffe's book \cite{GlimmJaffe} and the elementary properties of the Gohberg-Krein ${\det}_2$ determinant on the space of Hilbert-Schmidt operators imply that one has the equality of analytic functions
\begin{equation} \label{IdentityDet2}
\IE\Big[e^{-\lambda\langle\Wick{\phi_\eps^2},1\rangle}\Big] = {\det}_2\Big(\text{Id} + \lambda e^{2\eps\Delta}(H+c)^{-1}\Big)^{-1/2}
\end{equation}
on the disc $\{\vert\lambda\vert < \Vert e^{-2r\Delta}H^{-1}\Vert_{\text{HS}}\}\subset\IC$ with $\Vert\cdot\Vert_{\text{HS}}$ the Hilbert-Schmidt norm. For $\eps>0$ fixed, the analytic continuation property of the Gohberg-Krein determinant tells us that both sides of the equation extend as a meromorphic function over all of $\bbC$. We now prove that both terms converge to the correct limit to prove the resut.

\smallskip

We first take care of the probabilistic convergence of $\Wick{\phi^2_\eps}$ before looking at the partition function. For $p\ge2$ a large integer, we consider the convergence in $\CB_{2p,2p}^{-\delta}(\CS)$ for $\delta>0$ and conclude with the Besov embedding 
$$
\CB_{2p,2p}^{-\delta}(\CS) \hookrightarrow\CB_{\infty,\infty}^{-\delta-\frac{2}{p}}(\CS)
$$ 
in two dimensions. For $0<\eps_1,\eps_2\leq 1$, hypercontractivity ensures that
$$
\bbE\Big[\Vert \Wick{\phi^2_{\eps_1}} - \Wick{\phi^2_{\eps_2}}\Vert_{B^{-\delta}_{2p,2p}}^{2p}\Big] \lesssim \sum_{j\geq -1} 2^{-2pj\delta}\left(\int_\mcS \bbE\Big[P_j\big(\Wick{\phi^2_{\eps_1}} - :\Wick{\phi^2_{\eps_2}}\big)(x)^2\Big]\drm x \right)^p
$$
so it suffices to see that one has an $x$-uniform bound
\begin{equation} \label{EqUniformSmallBoundPhiSquare}
\bbE'\Big[P_j\big(\Wick{\phi^2_{\eps_1}} - \Wick{\phi^2_{\eps_2}}\big)(x)^2\Big] = o_{\eps_1,\eps_2}(1)
\end{equation}
as $\eps_1>0$ and $\eps_2>0$ go to $0$. Using the definition of Littlewood-Paley blocks from Appendix \ref{SubsectinoGeometricLP}, we get
\begin{equation*} \begin{split}
\bbE\Big[P_j&\big(\Wick{\phi^2_{\eps_1}} - \Wick{\phi^2_{\eps_2}}\big)(x)^2\Big] \\
&= \int_{\mcS\times\mcS} \bigg\{2\big(e^{\eps_1\Delta}(H+c)^{-1}e^{\eps_1\Delta}(z_1,z_2)\big)^2 + 2\big(e^{\eps_2\Delta}(H+c)^{-1}e^{\eps_2\Delta}(z_1,z_2)\big)^2   \\
&\qquad\qquad\qquad -2\big(e^{\eps_1\Delta}(H+c)^{-1}e^{\eps_2\Delta}(z_1,z_2)\big)^2 - 2\big(e^{\eps_2\Delta}(H+c)^{-1}e^{\eps_1\Delta}(z_1,z_2)\big)^2 \bigg\}   \\
&\qquad\qquad\times P_j(x,z_1)P_j(x,z_2) \drm z_1\drm z_2.
\end{split} \end{equation*}
We first start with the decomposition
\begin{eqnarray*}
e^{\eps_1\Delta}(H+c)^{-1}e^{\eps_2\Delta}(x,y)=e^{\eps_1\Delta}\left(\int_0^1e^{-t(H+c)} \drm t\right)e^{\eps_2\Delta} +  e^{\eps_1\Delta}\left(\int_1^\infty e^{-t(H+c)} \drm t \right)e^{\eps_2\Delta}.
\end{eqnarray*} 
Writing 
$$
\int_1^\infty e^{-t(H+c)} \drm t=e^{-\frac{1}{4}(H+c)}\left(\int_1^\infty e^{-(t-\frac{1}{2})(H+c)} \drm t\right) e^{-\frac{1}{4}(H+c)}
$$
with 
$$ 
e^{-(t-\frac{1}{2})(H+c)} : L^2(\mcS)\rightarrow L^2(\mcS)
$$ 
with operator norm bounded by $e^{-(t-\frac{1}{2})k}$ for $k>0$, we see that 
$$
\int_1^\infty e^{-(t-\frac{1}{2})(H+c)} \drm t = \mathcal{O}_{\mcB(L^2,L^2)}(1).
$$  
Since the operator $e^{-\frac{1}{4}(H+c)}$ has continuous positive kernel the map 
$$
x\in S\mapsto e^{-\frac{1}{4}(H+c)}(x,\cdot)\in L^2(\mcS)
$$ 
is continuous therefore we deduce that the composite operator
$$
e^{-\frac{1}{4}(H+c)}\left(\int_1^\infty e^{-(t-\frac{1}{2})(H+c)} \drm t\right) e^{-\frac{1}{4}(H+c)} 
$$
has a continuous Schwartz kernel. This means that one has the convergence
$$
e^{\eps_1\Delta}\left(\int_1^\infty e^{-t(H+c)} \drm t \right)e^{\eps_2\Delta} \underset{\eps_1,\eps_2\rightarrow 0}{\longrightarrow} \int_1^\infty e^{-t(H+c)} \drm t\in C^0(\mcS\times\mcS).
$$
Consider now the term $\int_0^1e^{-t(H+c)}\drm t$ which decomposes as
$$
\int_0^1e^{-t(H+c)}\drm t = \int_0^1 \left(e^{t(\Delta+c)}+A_te^{-tc}\right) \drm t.
$$
Since $A_t(x,y) = \mathcal{O}\big(t^{-1+1/q-\frac{\delta_0+\eta}{2}}\big)$ and $\frac{\delta_0+\eta}{2} <\frac{1}{q}$, the function $\int_0^1 A_te^{-tc} \drm t \in C^0(\mcS\times\mcS)$ converges with a continuous kernel and
$$
e^{\eps_1\Delta}\left(\int_0^1 A_te^{-tc} \drm t\right)e^{\eps_2\Delta} \underset{\eps_1,\eps_2\rightarrow 0}{\longrightarrow}\int_0^1 A_te^{-tc} \drm t \in C^0(\mcS\times\mcS). 
$$
It remains to observe that since the only `singular' term in
\begin{equation*} \begin{split}
B_{\eps_1,\eps_2}(z_1,z_2) &\defeq 2\big(e^{\eps_1\Delta}(H+c)^{-1}e^{\eps_1\Delta}(z_1,z_2)\big)^2 + 2\big(e^{\eps_2\Delta}(H+c)^{-1}e^{\eps_2\Delta}(z_1,z_2)\big)^2
\\
&\qquad -2\big(e^{\eps_1\Delta}(H+c)^{-1}e^{\eps_2\Delta}(z_1,z_2)\big)^2 - 2\big(e^{\eps_2\Delta}(H+c)^{-1}e^{\eps_1\Delta}(z_1,z_2)\big)^2
\end{split} \end{equation*}
is of the form $\int_0^1 e^{-(t+\eps_1+\eps_2)\Delta}(z_1,z_2) \drm t$, we have the convergence 
$$
\underset{\eps_1,\eps_2\rightarrow 0}{\lim}\,B_{\eps_1,\eps_2}(z_1,z_2)=0
$$ 
in $C^0(\mcS\times\mcS)$. We recall in identity \eqref{EqRepresentationLPProjector} of Appendix \ref{SubsectinoGeometricLP} that the kernels $P_j$ satisfy some identities of the form
\begin{eqnarray*}
P_j(x,y)=2^{j(\frac{d}{2}-1)}K_j(x,2^{\frac{j}{2}}(x-y))
\end{eqnarray*}
in well-chosen charts $U\times U$, where the kernels $K_j$ belong to a bounded family of smooth functions. It follows that one has
\begin{eqnarray*}
\left\vert\int_{U\times U} B_{\eps_1,\eps_2}(z_1,z_2)P_j(x,z_1)P_j(x,z_2) \drm z_1\drm z_2\right\vert \leq C 2^{-2j}\Vert B_{r_1,r_2} \Vert_{C^0(\mcS\times\mcS)}\underset{\eps_1,\eps_2\rightarrow 0}{\longrightarrow} 0  
\end{eqnarray*} 
where a positive constant $C$ independent of $j,\eps_1,\eps_2$. This concludes the proof of the bound \eqref{EqUniformSmallBoundPhiSquare}.

\smallskip

For the convergence of the partition function, define the joint variable
$$
\textbf{\textsf{X}}(\phi) \defeq \big(\phi,\Wick{\phi^2}\big)\in H^{-\delta}(\mcS)\times H^{-2\delta}(\mcS)
$$
for $\delta>0$ and equip the product space $H^{-\delta}(\mcS)\times H^{-2\delta}(\mcS)$ with the metric
$$
\llparenthesis(a,b)\rrparenthesis \defeq \Vert a\Vert_{H^{-\delta}} + \Vert b\Vert_{H^{-2\delta}}^{1/2}.
$$
We consider $\textbf{\textsf{X}}$ as a measurable function of $\phi$. The Cameron-Martin embedding $\CC\CM\subset H^{1-\kappa}(\CS)$ implies that almost surely one has for all $h\in\CC\CM$
$$
\textbf{\textsf{X}}(\phi+h) = \textbf{\textsf{X}}(\phi) + 2h\phi + h^2,
$$
with a well-defined product $h\phi$. The function $\llparenthesis\textbf{\textsf{X}}(\cdot)\rrparenthesis$ satisfies then $\phi$-almost surely the estimate
\begin{equation} \label{EqInequalityFernique}
\llparenthesis\textbf{\textsf{X}}(\phi)\rrparenthesis \lesssim \llparenthesis\textbf{\textsf{X}}(\phi-h)\rrparenthesis + \Vert h\Vert_{\CC\CM}
\end{equation}
for all $h\in\CC\CM$ for an absolute implicit multiplicative constant in the inequality. One then gets from Friz and Oberhauser generalized Fernique's theorem \cite{FrizOberhauser} that the random variable $\llparenthesis\textbf{\textsf{X}}(\phi)\rrparenthesis$ has a Gaussian tail. The random variable $\exp\left(-\lambda\langle\Wick{\phi^2},1\rangle\right)$ is thus integrable for $\lambda\in\bbC$ small enough. If one defines similarly
$$
\textbf{\textsf{X}}_\eps(\phi) \defeq \big(\phi_\eps,\Wick{\phi_\eps^2}\big)\in H^{-\delta}(\mcS)\times H^{-2\delta}(\mcS),
$$
then the function $\llparenthesis\textbf{\textsf{X}}_\eps(\cdot)\rrparenthesis$ also satisfies the estimate
$$
\llparenthesis\textbf{\textsf{X}}_\eps(\phi)\rrparenthesis \lesssim \llparenthesis\textbf{\textsf{X}}_\eps(\phi-h)\rrparenthesis + \Vert h\Vert_{\CC\CM}
$$
with the same implicit constant as in \eqref{EqInequalityFernique}. The conclusion of Fernique's generalized theorem is quantitative and can be written in terms of the $\overline{\textrm{erf}}$ function
$$
\overline{\textrm{erf}}(z) = 1-\textrm{erf}(z) = \frac{1}{\sqrt{2\pi}}\int_z^\infty e^{-a^2/2}da.
$$
If one sets 
$$
\mu_{a,\eps} \defeq \bbP'\big(\llparenthesis\textbf{\textsf{X}}_\eps(\phi)\rrparenthesis \leq a\big), \qquad a'_\eps \defeq \textrm{erf}^{-1}\,(\mu_{a,\eps}),
$$
for a fixed $a>0$ such that $0<\mu_{a,\eps}<1$, then 
$$
\bbP'\big(\llparenthesis\textbf{\textsf{X}}_\eps(\phi)\rrparenthesis > m \big), \leq \overline{\textrm{erf}}(a_\eps'+\sigma m),
$$
for a positive constant $\sigma$ that depends only on $a$ and the implicit constant in \eqref{EqInequalityFernique}. As $\llparenthesis\textbf{\textsf{X}}_\eps(\cdot)\rrparenthesis$ is converging to $\llparenthesis\textbf{\textsf{X}}(\cdot)\rrparenthesis$ in $L^2$ as a random variable conditionned to $\xi$, one can choose a constant $a$ such that $\bbP\big(\llparenthesis\textbf{\textsf{X}}(\cdot)\rrparenthesis \leq a\big)$ is also in $(0,1)$. It is thus possible to find an $a'$ such that one has 
$$
\sup_{0<\eps\leq 1}\,\bbP\big(\llparenthesis\textbf{\textsf{X}}_\eps(\phi)\rrparenthesis > m \big) \leq \overline{\textrm{erf}}(a'+\sigma m).
$$
It follows from that estimate that the family of random variables $\exp\left(-\lambda\langle\Wick{\phi_\eps^2},1\rangle\right)$ for $0<\eps\leq 1$ and $\lambda$ in a small ball of $\bbC$, is uniformly integrable; so it converges in $L^1(\Omega',\bbE')$ to $\exp\hspace{-0.05cm}\left(-\lambda\langle\Wick{\phi^2},1\rangle\right)$. 

\smallskip

For the convergence of the determinant, we have that the operators $(H+c)^{-1}e^{-\eps\Delta}(H+c)^{-1}$ are indeed trace class as symmetric non-negative operators with kernels $K_\eps(x,y)$ satisfying the estimate
$$
\int_\mcS K_\eps(x,x)\,\mu(\drm x) < \infty
$$
uniformly in $\eps\in[0,1]$ using the estimate on the Green function $G$ from Proposition \ref{PropHeatKernalBounds}. It follows that 
\begin{equation*} \begin{split}
\text{Tr}_{L^2}\Big(&(H+c)^{-1}\big(e^{2\eps\Delta}-1\big)\big(e^{2\eps\Delta}-1\big)(H+c)^{-1}\Big)   \\
&= \text{Tr}_{L^2}\Big((H+c)^{-1}e^{4\eps\Delta}H^{-1}\Big) - 2\,\text{Tr}_{L^2}\Big(H^{-1}e^{2\eps\Delta}(H+c)^{-1}\Big) + \text{Tr}_{L^2}\big((H+c)^{-2}\big)   \\
&= \int_\mcS G(x,y)p^\Delta_{4\eps}(y,z)G(z,x)\drm z\drm y\drm x - 2\int_\mcS G(x,y)p^\Delta_{2\eps}(y,z)G(z,x)\drm z\drm y\drm x + \int_\mcS G(x,y)^2\drm x
\end{split} \end{equation*}
is converging to $0$. The continuity of the ${\det}_2$ function on the ideal of Hilbert-Schmidt operators on $L^2(\mcS)$ then gives the convergence of the determinant hence
\begin{equation}
\IE\Big[e^{-\lambda\langle\Wick{\phi^2},1\rangle}\Big] = {\det}_2\Big(\text{Id} + \lambda (H+c)^{-1}\Big)^{-1/2}
\end{equation}
on the disc $\{\vert\lambda\vert < \Vert H^{-1}\Vert_{\text{HS}}\}\subset\IC$. Since the analytic continuation to all of $\bbC$ of the locally defined function $\lambda\mapsto{\det}_2\big(\textrm{Id} + \lambda (H+c)^{-1}\big)$ has its zero set equal to 
$$
\big\{-z^{-1}\,;\,z\in\sigma((H+c)^{-1})\big\},
$$
we see that the partition function $Z(\cdot)$ determines the spectrum of $H+c$, hence the spectrum of $H$. The formula involving the $a_n$ is obtain with the general identity
$$
{\det}_2(1+\lambda A) = \exp\Big(-\sum_{n\geq 2} \frac{(-\lambda)^n}{n}\,\text{Tr}(A^n)\Big)
$$
which follows from the fact that $\det(e^B)=e^{\text{TR(B)}}$ with $B=\log(1+\lambda A)$ as a power serie valid for any Hilbert-Schmidt operator $A$ on $L^2(\mcS)$, see again Chapter 9 in Glimm and Jaffe's book \cite{GlimmJaffe}.
\end{proof}

\medskip

The proof of Theorem \ref{ThmWickCharacterizationSpectrumH} actually tells us that for every non-negative function $f$ in $B^{1/p}_{p,\infty}(\mcS)$ with  $1/p>2\nu$, one has the formula
\begin{equation} \label{EqPartitionFunctionf} \begin{split}
Z(f) \defeq \bbE'\left[e^{-:\phi^2:(f)}\right] = {\det}_2\Big(\textrm{Id} + M_{f^{1/2}} (H+c)^{-1}M_{f^{1/2}}\Big)^{-1/2}.
\end{split} \end{equation}
Indicators of subsets of $\mcS$ with finite perimeter are elements of the spaces $B^{1/p}_{p,\infty}(\mcS)$ with $1/p>2\nu$ -- see e.g. Theorem 2 in Sickel's survey \cite{Sickel}. 

To emphasize that the real-valued quantities $Z(\lambda)$ and $a_n$ are random and their laws depend on the Riemannian metric space $(\mcS,g)$ we write $Z(\lambda)(\mcS,g)$ and $a_n(\mcS,g)$. The next statement gives a characterization of the law of the spectrum of $H$, a function of $(\mcS,g)$, in terms of the law of the $a_n(\mcS,g)$. Write here $H(\mcS,g)$ to emphasize this dependence.

\medskip

\begin{corollary} \label{CorSpectralLaw}
Let $(\mcS_1,g_1)$ and $(\mcS_2,g_2)$ be two Riemannian closed surfaces. Then the spectra of the operators $H(\mcS_1,g_1)$ and $H(\mcS_2,g_2)$ have the same law iff the sequences $\big(a_n(\mcS_1,g_1)\big)_{n\geq 2}$ and $\big(a_n(\mcS_2,g_2)\big)_{n\geq 2}$ have the same law.
\end{corollary}

\medskip

Either condition is equivalent to the fact that the functions $Z(\cdot)(\mcS_1,g_1)$ and  $Z(\cdot)(\mcS_2,g_2)$ have the same law.

\medskip

\begin{proof}
Use Skorohod representation theorem to turn equality in law into almost sure equality on a different probability space.

If the two sequences $\big(c_n(\mcS_1,g_1)\big)_{n\geq 2}$ and $\big(c_n(\mcS_2,g_2)\big)_{n\geq 2}$ are equal the two functions $Z(\cdot)(\mcS_1,g_1)$ and  $Z(\cdot)(\mcS_2,g_2)$ are equal, and the functions ${\det}_2\big(1+\lambda H(\mcS_1,g_1)\big)$ and ${\det}_2\big(1+\lambda H(\mcS_2,g_2)\big)$ of $\lambda$ coincide on a small disk, hence on all of $\bbC$. Given the relation between the zero set of these functions and the spectrum of the operators $H(\mcS_1,g_1)$ and $H(\mcS_2,g_2)$ these spectra need to coincide. The function $Z$ is determined by the spectrum of $H$ since the $a_n$ has that property from \eqref{EqTraceFormulaAn}.
\end{proof}

\medskip

Corollary \ref{CorSpectralLaw} somehow says that the law of the partition function of $:\hspace{-0.07cm}\phi^2\hspace{-0.07cm}:$ determines the law of the spectrum of $H$.

\medskip

\begin{remark} 
The Anderson Gaussian free field introduced in this section is a new object. It echoes some other works that somewhat share a similar spirit. In Caravenna, Sun and Zygouras' work \cite{CSZ16} and Bowditch and Sun's work \cite{BowditchSun} the authors consider the scaling limit of an Ising model on $\mathbb{Z}^2$ at the critical temperature subject to some random singular magnetic field modelled by white noise. From a constructive quantum field theory viewpoint this is similar to studying some $\phi^4_2$ measure with source term
\begin{eqnarray*}
\mathbb{E}\left( e^{-\int_{\mcS}\left( :\phi^4:-\lambda_c:\phi^2:\right)\mu +\int_{\mcS} \phi\xi\mu } \right)
\end{eqnarray*}
where $\lambda_c\geqslant 0$ is chosen to be the critical parameter of the $\phi^4_2$ measure -- it plays the role of the critical temperature in the Ising model, with a white noise source term $\xi$, and where the expectation is taken with respect to a particular massive Gaussian free field measure. The existence of the critical value $\lambda_c$ follows from the work of Glimm, Jaffe and Spencer \cite{GlimmJaffe}. In our case, we study a free field where white noise plays the role of a random singular potential instead of a random magnetic field.
\end{remark}

\section{Polymer measure and Anderson diffusion}
\label{SectionPolymer}

The semigroup $e^{-tH}$ is not conservative hence the Anderson heat kernel $p_t(x,\cdot)$ is not of unit mass over $\CS$. Still there are two natural processes that one can consider given such a Schrödinger operator, namely its associated polymer measure and its ground state diffusion.

\subsection{Anderson polymer measure}

Given a smooth potential $V$, the polymer measure on $[0,T]$ of length $T>0$ starting at $x\in\CS$ is given by the measure
\begin{equation}
\CQ_x^T(\drm X)=\frac{1}{Z_T}e^{\int_0^TV(X_t)\drm t}\CP_x^T(\drm X)
\end{equation}
with $\CP_x^T$ the law of the Brownian motion on $\CS$ starting at $x$ stopped at time $T$. The Wiener measure $\CP_x^T$ is penalized according to the potential $V$ and the path has a higher probability of lying where the potential is low; the constant $Z_T>0$ is a normalization constant that depends on the length $T$ to ensure that $\CQ_x^T$ is a probability measure. In the context of a rough potential, it is unclear how to interpret the evaluation of the potential $V(X_t)$ and then that the normalization constant $Z_T$ is finite. The Feynman-Kac formula relates such polymer measure to the associated stochastic heat equation
\begin{equation}
\partial_tu=\Delta u+Vu,
\end{equation}
and this was the starting point of the construction by Alberts, Khanin and Quastel \cite{AKQ} for the polymer measure in one dimension with a spacetime Gaussian white noise. We follow their approach since the Anderson heat kernel precisely gives the probability transition of the underlying process. Cannizzaro and Chouk \cite{CannizzaroChouk} constructed the polymer measure on the two dimensional torus using the KPZ equation with a Girsanov transform relying on SDEs with time dependent drift. In relation with the Anderson diffusion, one could interpret our construction of the polymer as a relation to an SDE with a time independent distributional drift. See also the recent work by Berger and Lacoin \cite{BergerLacoin} for a different approach to construct the polymer measure in a random environment given by a Lévy noise.

\medskip

\begin{definition}
The \textbf{polymer measure} is the measure with finite dimensional projections 
\begin{equation}
\CQ_x^T(X_{t_1}\in A_1,\ldots,X_{t_n}\in A_n)=\frac{1}{p_T(x)}\int_{A_1\times\ldots\times A_n\times\CS}\prod_{i=0}^np_{t_{i+1}-t_i}(x_i,x_{i+1})\drm x_{i+1}
\end{equation}
for any $t_0=0<t_1<\ldots<t_n<T=t_{n+1}$ and $A_1,\ldots,A_n\subset\CS$ measurable sets with $x_0=x$, and normalization constant 
\begin{equation}
p_T(x)=\int_\CS p_T(x,y)\drm y
\end{equation}
depending on the initial point $x\in\CS$ and the length $T>0$ of the polymer.
\end{definition}

\medskip

The measure $\CQ_x^T$ denotes the law of the random continuum polymer fixed at $x\in\CS$ in $t=0$ of length $T>0$. (One could also consider the measure $\CQ_{x,y}^T$ with the additional constraint of being fixed at $y\in\CS$ in $t=T$.) The endpoint of our polymer is free thus the final integration over $\CS$ in the previous definition. While the measure $\CQ_x^T$ is uniquely characterized by its finite dimensional projection, it is not clear a priori if the previous measure is supported on the set of continuous paths $C([0,T],\CS)$. This is granted by the Gaussian upper bounds that we have on the Anderson heat kernel, together with Kolmogorov Theorem.

\smallskip

\medskip

Since the polymer measure is formally given by the expression
\begin{equation}
\CQ_x^T(\drm X)=\frac{1}{Z_T}e^{\int_0^TV(X_t)\drm t}\CP_x^T(\drm X)
\end{equation}
for any $T>0$ and $x\in\CS$, it is natural to ask if it is absolutely continuous with respect to the Wiener measure $\CP_x^T$. This is not the case and one can prove that the measures are singular. This was already proved by Cannizzaro and Chouk, see Theorem 1.4 \cite{CannizzaroChouk}. We give a proof in our context for completeness.

\medskip

\begin{proposition} \label{PropSingularLaws}
For any $x\in \mcS$ and $T>0$, the polymer measure $\CQ_x^T$ is $\bbP$-almost surely singular with respect to the Wiener measure $\CP_x^T$.
\end{proposition}

\medskip

\begin{proof}
Let $(\eps_n)_{n\ge0}\subset(0,1]$ be a sequence decreasing to $0$ and consider
\begin{equation}
D_n(B):=e^{-\int_0^T(\xi_\eps+\frac{\log\eps}{4\pi})(B_t)\drm t}
\end{equation}
which corresponds to the density of the polymer with regularized noise with respect to the Wiener measure. We prove that the event
\begin{equation}
\limsup_n\{D_n<1\}
\end{equation}
is of probability $1$ for $\CP_x^T$ and $0$ for $\CQ_x^T$.

\smallskip

For the first part, Feynman-Kac formula for smooth potential allows to compute the expectation of such quantities, we have
\begin{equation}
\IE_x\big[D_n^{\frac{1}{2}}\big]=\big(e^{-T(\Delta+\frac{1}{2}\xi_{\eps_n}+\frac{\log\eps_n}{8\pi})}1\big)(x)
\end{equation}
where $\IE_x$ denotes expectations with respect to $\CP_x^T$. Since the renormalization constant is quadratic with the noise, we write
\begin{equation}
\IE_x\big[D_n^{\frac{1}{2}}\big]=e^{-T\frac{\log\eps_n}{16\pi}}\big(e^{-T(\Delta+\frac{1}{2}\xi_{\eps_n}+\frac{\log\eps_n}{16\pi})}1\big)(x)
\end{equation}
which converges to $0$ as $n$ goes to infinite since the semigroup converges to a finite quantities. So
$$
\lim_{n\to\infty} \bbE_x\Big[D_{\eps_n}^{1/2}\Big] = 0
$$
and $\CP_x^T(D_{\eps_n} > 1)$ tends to $0$ from Chebychev inequality. One has as a consequence
$$
\CP_x^T\Big(\limsup_n \big\{D_{\eps_n} < 1\big\}\Big) \ge \limsup_n \CP_x^T(D_{\eps_n} < 1) = 1.
$$
For the second part, we have
$$
\CQ_x^T(D_{\eps_k}<1)\le\liminf_n\CQ_{\eps_n,x}^T(D_{\eps_k}<1),
$$
for fixed $k\ge0$ with $\CQ_{\eps_n,x}^T$ the polymer measure with the potentiel $\xi_{\eps_n}$. Using that $\CD_{\eps_n,x}^T$ has a density with respect to $\CP_x^T$, we get
\begin{align*}
\CQ_{\eps_n,x}^T(D_{\eps_k}<1) &= \bbE_x\Big[e^{-\int_0^T (\xi_{\eps_n}+(\log \eps_n)/(4\pi))(B_t)\drm t}\IDC_{D_{\eps_k}<1}\Big]\\
&\le\bbE_x\Big[e^{-\int_0^T (\xi_{\eps_n}+(\log \eps_n)/(4\pi))(B_t)\drm t}D_{\eps_k}^{-1/2}\IDC_{D_{\eps_k}<1}\Big]\\
&\leq  \bbE_x\Big[e^{-\int_0^T[\xi_{\eps_n}+(\log \eps_n)/(4\pi) - 1/2(\xi_{\eps_k}+(\log \eps_k)/(4\pi))](B_t)\drm t}\Big]   \\
&\le e^{-T(\log \eps_k)/(16\pi)}\, \bbE_x\Big[e^{-\int_0^T[\xi_{\eps_n} +(\log \eps_n)/(4\pi) - (1/2\xi_{\eps_k}+(\log \eps_k)/(16\pi))](B_t)\drm t}\Big].
\end{align*}
As 
$$
{\sf \Pi}\Big(X_{\eps_n}+\frac{1}{2}X_{\eps_k}, \xi_{\eps_n}+\frac{1}{2}\xi_{\eps_k}\Big) - \frac{\log \eps_n}{4\pi} + \frac{5}{4}\,\frac{\log \eps_k}{4\pi}
$$
converges in probability in $C^{1-2\kappa}(\mcS)$ as $n$ goes to $\infty$ then $k$ goes to $\infty$, one sees that the quantity
$$
\bbE_x\Big[e^{-\int_0^T[\xi_{\eps_n} +(\log \eps_n)/(4\pi) - (1/2\xi_{\eps_k}+(\log \eps_k)/(16\pi))-(\log \eps_k)/(4\pi)](B_t)\drm t}\Big]
$$
converges as $n$ goes first to $\infty$ then $k$ goes to $\infty$. It follows that 
$$
\CQ_{\eps_n,x}^T(D_{\eps_k}<1) \lesssim e^{\frac{3T}{4}\frac{\log \eps_k}{4\pi}}
$$
uniformly in $n$ and $k$ so
$$
\CQ_x^T(D_{\eps_k}<1) \lesssim e^{\frac{3T}{4}\,\frac{\log \eps_k}{4\pi}}.
$$
Choosing a sequence $\eps_k$ that decreases sufficiently fast to $0$ provides then an upper bound for $\CQ_{\eps_n,x}^T(D_{\eps_k}<1)$ that allows to conclude with Borel-Cantelli lemma that 
$$
\CQ_x^T\Big(\limsup_k\{D_{\eps_k}<1\}\Big) = 0.
$$
\end{proof}

\subsection{Anderson diffusion}

Another natural process associated to a Schrödinger operator with a spectral gap $\lambda_0<\lambda_1$ and a ground state $\Psi>0$ is via the infinitesimal generator
\begin{equation}
\SL=\Psi^{-1}(H-\lambda_0)\Psi
\end{equation}
with corresponding heat semigroup 
\begin{equation}
e^{-t\SL}=e^{t\lambda_0}\Psi^{-1}e^{-tH}\Psi
\end{equation}
for any $t\ge0$. It has a kernel $p_t^\SL(x,y)$ given by
\begin{equation}
p_t^\SL(x,y) = e^{t\lambda_0} \, \frac{\Psi(y)}{\Psi(x)} \, p_t(x,y)
\end{equation}
given the heat kernel $p_t$ associated to the Schrödinger operator. This is indeed a conservative semigroup
\begin{align}
\int_\CS p_t^\SL(x,y)\mu(\drm y)&=\frac{e^{t\lambda_0}}{\Psi(x)}\int_\CS p_t(x,y)\Psi(y)\mu(\drm y) = \frac{e^{t\lambda_0}}{\Psi(x)}(e^{-tH}\Psi)(x) = 1,
\end{align}
since $H\Psi=\lambda_0\Psi$. The operator $\SL$ already appeared in the proof of the two-sided Gaussian bounds on $p_t$, in Proposition \ref{PropHeatKernalBounds}, with
\begin{equation}
\SL=-\Delta-2\nabla\log(\Psi)\nabla
\end{equation}
in the case of the Anderson ground state $\Psi=v_0$. The operator $-\frac{1}{2}\SL$ corresponds to the infinitesimal generator of the SDE
\begin{equation}
\drm X_t=(\nabla\log\Psi)(X_t)\drm t+\drm B_t
\end{equation}
which is an SDE with time independent distributional drift for the Anderson diffusion since $u_0\in C^{1-\kappa}(\CS)$.

\medskip

\begin{definition}
The \textbf{Anderson diffusion} is the process defined by the conservative semigroup $(e^{-\frac{1}{2}\SL})_{t>0}$ generated by the operator $-\frac{1}{2}\SL$.
\end{definition}

\medskip

In comparison to the Anderson polymer, the finite dimensional projection of the Anderson diffusion starting at $x\in\CS£ $ are given by
\begin{align}
\CP_x(X_{t_1}\in A_1,\ldots,X_{t_n}\in A_n)&=\int_{A_1\times\ldots\times A_n}\prod_{i=0}^{n-1}p_{t_{i+1}-t_i}^\SL(x_i,x_{i+1}) \, \drm x_{i+1}   \\
&=\int_{A_1\times\ldots\times A_n\times\CS}e^{-T\lambda_0} \, \frac{\Psi(x_n)}{\Psi(x)} \, \prod_{i=0}^{n-1}p_{t_{i+1}-t_i}(x_i,x_{i+1}) \, \drm x_{i+1}
\end{align}
for any $t_0=0<t_1<\ldots<t_n<T=t_{n+1}$ and $A_1,\ldots,A_n\subset\CS$ measurable sets with $x_0=x$. Note that the Anderson diffusion stopped at time $T>0$ gives a different measure than the polymer measure as it is a Markov process while the polymer is not. Like for the polymer measure, the upper bound on the Anderson heat kernel immediately gives that this is indeed a measure supported on continuous paths in $\CS$. Our analysis of the Anderson heat kernel gives a number of properties of the two paths measures.

\medskip

\begin{proposition} \label{PropQuadraticVariationProcess}
The Anderson polymer and the Anderson diffusion are almost surely $\alpha$-Hölder for any $\alpha<\frac{1}{2}$ as paths with values in $\CS$. Moreover the paths are almost surely of finite quadratic variation.
\end{proposition}

\medskip

\begin{proof}
The Hölder regularity follows from Kolmogorov criterion, the Gaussian upper bounds on the Anderson heat kernel indeed gives that the paths are of the same regularity of the Brownian motion in $\CS$.

We now prove that the quadratic variation of the canonical process on path space is a well-defined random variable under $Q_x^T$. This means that 
$$
\sum_{i=0}^n d(w_{t_{i+1}},w_{t_i})^2
$$
converges in $L^2(Q_x^T)$ to (the constant random variable) $t$, for each $t$ when the mesh of a partition $0<t_1<\cdots<t_n<t$ of an interval $[0,t]$, with $t_0\defeq0$ and $t_{n+1}\defeq1$, goes to $0$. (Do not mingle the fact for a process to have a finite quadratic variation process and the property of its sample paths to be almost surely of finite $2$-variation. Brownian motion has for instance a finite quadratic variation process on any finite interval but has almost surely an infinite $2$-variation on any finite interval.) To prove the preceding convergence we notice that the fine asymptotic from Corollary \ref{CorAsymptoticKernel} gives
\begin{equation} \label{EqExpectedScaling}
\bbE_x\big[d(w_{t_{i+1}},w_{t_i})^2\big] = t_{i+1}-t_i + O(t_{i+1}-t_i)^b
\end{equation}
for a constant $b>1$, and that   
$$
\bbE_x\big[d(w_{t_{i+1}},w_{t_i})^4\big] = O(t_{i+1}-t_i)^b,
$$
from the Gaussian upper bound on the heat kernel. Chebychev's inequality then gives the result. We note here for later purposes that for each $t$, there is a sequence of partitions of the interval $[0,t]$ such that the corresponding sum of squared increments converges almost surely to $t$. The quadratic variation process thus depends only on the equivalence class of a finite non-negative measure on path space under the equivalence relation given by reciprocal absolute continuity.

Note that the Gaussian lower and upper estimates on the heat kernel $p_t$ proved in Proposition \ref{PropHeatKernalBounds} are not sufficient to get back the exact scaling relation \eqref{EqExpectedScaling}. One really needs the result of item {\it (1)} Theorem \ref{ThmAsymptotic} for that purpose.
\end{proof}

\subsection{Wick square of Anderson Gaussian free field and the Anderson diffusion}
\label{SectionDynkin}

The study of the links between some Markov fields and some Poissonian ensembles of Markov loops goes back to Symanzik' seminal work \cite{Symanzik}. It was elaborated in a large number of works and we take advantage here of the general result proved by Le Jan in \cite{LeJan}, giving a correspondence between the occupation measure of a loop ensemble and Wick square of some Gaussian free field -- see Section 9 therein. It allows at no cost to relate (a measure built from) the Anderson diffusion to the Wick square of the Anderson free field that was the object of Theorem \ref{ThmWickCharacterizationSpectrumH}. We dress the table before bringing the dishes.

\medskip

Rather than working with the polymer measure built from the operator $H-\lambda_0(\widehat{\xi}\,)$ we pick a positive constant $a$ and work with the operator built from $H-\lambda_0(\widehat{\xi}\,)+a$. With the notations of Section \ref{SectionAndersonGFF} one takes here $c=-\lambda_0(\widehat{\xi}\,)+a$. This choice ensures that the Green function of the corresponding semigroup is finite and has the properties stated and used in Section \ref{SectionAndersonGFF}. This amounts to adding killing at a constant rate $a$ for the Anderson diffusion. This does not change its properties and we have in particular that the corresponding diffusion paths have an associated quadratic variation process equal to the travelling time and defined on a random lifetime interval $[0,\zeta)$. Set 
$$
\frac{e^{t(\lambda_0(\widehat\xi\,)-a)} \, p_t(x,y) \, u_0(y)}{u_0(x)}
$$
and denote by $\overline{P}^t_{x,x}$ the unnormalized excursion measure of duration $t$ started from $x\in\mcS$. It is characterized by the identity
\begin{equation*} \begin{split}
\overline{P}^t_{x,x}\Big(X_{t_1}\in \mu(dx_1),\dots X_{t_k}\in \mu(dx_k)\Big) &= \overline{p}_{t_1}(x,x_1) \, \overline{p}_{t_2-t_1}(x_1,x_2)\dots \overline{p}_{t-t_k}(x_k,x) \, \mu(dx_1)\dots \mu(dx_k)   \\
&= p_{t_1}(x,x_1)\,p_{t_2-t_1}(x_1,x_2)\dots p_{t-t_k}(x_k,x) \, \mu(dx_1)\dots \mu(dx_k)
\end{split} \end{equation*}
for all $0\leq t_1\leq \cdots\leq t_k\leq t$. Note that these quantities are independent of $u_0$. This non-negative measure has a finite mass equal to $\overline{p}_t(x,x)$. A standard argument using the symmetry of $p_t(x,y)$ as a function $(x,y)$ shows that the measure $\overline{P}$ is supported on (rooted) loops of H\"older regularity strictly less than $1/2$. The loop measure is defined as
$$
\mathscr{M}(\cdot) \defeq \int_\mcS\int_0^\infty \frac{1}{t}\,\overline{P}^t_{x,x}(\cdot)\,dt\,\mu(dx).
$$
It follows from Proposition \ref{PropQuadraticVariationProcess} that the factor $1/t$ in this integral accounts for the intrinsic lifetime of the loop, the quadratic variation process -- so this non-negative measure is indeed a measure on unrooted loops. Note that it has an infinite mass that comes from the mass of small loops. Denote by $\bbE_{\mathscr{M}}$ the expectation operator associated with $\mathscr{M}$ and by $\zeta(\ell)$ the lifetime of a loop $\ell$. For such a loop we define a measure on $\mcS$ setting
$$
\widehat{\ell}(\cdot) \defeq \int_0^{\zeta(\ell)}\delta_{\ell(s)}(\cdot)\,ds.
$$
One has for any non-negative function $f$ on $\mcS$ and all $n\geq 1$
\begin{equation} \label{EqHatEll}
\bbE_{\mathscr{M}}\big[\widehat{\ell}(f)^n\big] = (n-1)!\,\int_{\mcS^n} G(x_1,x_2)f(x_2)G(x_2,x_3)f(x_3)\cdots G(x_n,x_1)f(x_1)\,\mu(dx_1)\dots \mu(dx_n),
\end{equation}
and 
\begin{equation} \label{EqLevyMeasure}
\bbE_{\mathscr{M}}\left[e^{-z\widehat{\ell}(f)} + z\widehat{\ell}(f) -1 \right] = -\log\textrm{det}_2\Big(\textrm{Id} + zM_{f^{1/2}}GM_{f^{1/2}}\Big),
\end{equation}
from an elementary series expansion and the preceding equality. We used here the same notation for the Green kernel $G$ of $H+c$ and its associated operator $(H+c)^{-1}$. Le Jan's proof \cite{LeJan} of identity \eqref{EqHatEll} applies verbatim here. The quantity that naturally appears in formula \eqref{EqHatEll} involves the Green function of the operator $u_0^{-1}e^{-t(H+c)}(u_0\cdot)$, that is the conjugate of $(H+c)^{-1}$ by the multiplication operator by $u_0$. The expression \eqref{EqHatEll} being cyclic in $(x_1,\dots,x_n)$ it turns out to be independent of $u_0$.

Given $\gamma\geq 0$ denote by $\Lambda_\gamma$ a Poisson process on the space of (unrooted) loops over $\mcS$ with intensity $\gamma\mathscr{M}$. It is characterized by its characteristic function
$$
\bbE\big[e^{i\Lambda_\gamma(F)}\big] = \exp\left(\gamma\int\big(e^{iF(\ell)}-1\big)\,\mathscr{M}(d\ell)\right),
$$
for all functions $F$ on loop space that are null on loops of sufficiently small lifetime -- so the resulting quantity $\Lambda_\gamma(F)$ is almost surely well-defined. Denote by $A_\gamma$ the support of $\Lambda_\gamma$, so $\Lambda_\gamma = \sum_{\ell\in A_\gamma} \delta_\ell$. The regularized renormalized occupation measure of $\Lambda_\gamma$ is defined for each $r>0$ as the non-negative measure on $\mcS$
$$
\mathcal{O}_\gamma^R(f) \defeq \sum_{\ell\in A_\gamma} \left( {\bf 1}_{\zeta(\ell)>r} \, \widehat{\ell}(f) - \gamma \, \bbE_{\mathscr{M}}\big[{\bf 1}_{\zeta(\ell')>r} \, \widehat{\ell'}(f)\big]\right);
$$
the expectation is over $\ell'$ and $f$ is a generic non-negative continuous function on $\mcS$. For $\gamma$ and $f$ fixed the continuous time random process $\gamma\mapsto \mathcal{O}_\gamma^R(f)$ is actually a L\'evy process with positive jumps with characteristic function
$$
\bbE\Big[e^{-\mathcal{O}_\gamma^R(f)}\Big] = \exp\left(-\gamma\,\bbE_{\mathscr{M}}\Big[ {\bf 1}_{\zeta(\ell')>r}\big(e^{-\widehat{\ell}(f)} + \widehat{\ell}(f) -1\big)\Big]\right)
$$
converging to its natural limit as $r$ goes to $0$. The limit L\'evy process is denoted by $(\mathcal{O}_\gamma(f))_{\gamma\geq 0}$. (All this is explained in detail in Le Jan's work \cite{LeJan}.) The following result follows from the preceding analysis and the formulae \eqref{EqLevyMeasure} and \eqref{EqPartitionFunctionf} for the partition function of the Wick square of the Anderson Gaussian free field.

\medskip

\begin{theorem}
For every continuous function $f$ on $\mcS$ that is also in $B^{1/p}_{p,\infty}(\mcS)$, with $1/p<2\nu$, one has the identity
$$
\bbE\big[e^{-\mathcal{O}_{1/2}(f)}\big] = \bbE\big[e^{-:\phi^2\hspace{-0.03cm}:\hspace{-0cm}(f)}\big].
$$
\end{theorem}

One deduces from this identity that the renormalized occupation measure of the loop measure of polymer paths has the same distribution as the Wick square of the Anderson Gaussian free field. It has in particular a version that has almost surely regularity $-2\nu$ in the Sobolev scale. This identification does not tell us that $\mathcal{O}_{1/2}$ is a measure, despite its name.

\medskip

\appendix

\section{Meromorphic Fredholm theory with a parameter}
\label{SectionFredholm}

We prove Theorem \ref{ThmFredholm} in this section. As a guide to the subject of this appendix, the reader will find in Appendix D of Zworski's book \cite{Zworski} an elementary account of the usual, parameter-free, meromorphic Fredholm theory.

\medskip

\begin{proof}
Our proof follows closely the proof given by Borthwick in Theorem 6.1 of \cite{Borthwick}. It suffices to prove the result near any $z_0\in U$ which contains only finitely many poles of $K$. With this assumption, we may decompose 
$$
K(z,e) = A(z,e) + F(z,e),
$$ 
where $F(z,e)$ is a meromorphic family of finite-rank operators for $z\in U$ and $A(z,e)$ is a holomorphic family of compact operators. Both operators depend continuously on the parameter $e$. Using the approximation of the compact operator $A(z_0,e)$ by finite-rank operators, and assuming $U$ is sufficiently small and that we choose a sufficiently small neighbourhood of $e_0$, we can find a fixed finite-rank operator $B$ such that
$$
\big\Vert A(z,e) - B \big\Vert < 1
$$
for all $z\in U$. Note that  implies that $\textrm{Id}-A(z,e) + B$ is holomorphically invertible for $z\in U$, by the usual Neumann series as
$$
\big(\textrm{Id}-A(z,e) + B \big)^{-1}=\sum_{k=1}^\infty \big(A(z,e) - B\big)^k.
$$
Since the Neumann series converges absolutely in $\mathcal{B}(\mathcal{H},\mathcal{H})$ uniformly in $(z,e)$ in some neighborhood of $(z_0,e_0)$ and each term $\left(A(z,e) - B\right)^k$ is continuous in $u$, it follows that the map 
$$
e\mapsto  \left(\textrm{Id}-A(z,e) + B\right)^{-1}\in \mathcal{B}(\mathcal{H},\mathcal{H})
$$ 
is continuous. Thus if we set
$$
G(z,e) \defeq \big(F(z,e) + B\big)\,\big(\textrm{Id}-K(z,e) + B\big)^{-1} 
$$
then we can write
$$
\textrm{Id}-K(z,e) = \big(\textrm{Id}- G(z,e)\big) \, \big(\textrm{Id}-K(z,e) + B\big)^{-1}.
$$
It is immediate that $G(z,e)$ has finite rank and depends continuously on $e$ by its construction involving the finite rank operators $F(z,e), B$. We already know that $\left(\textrm{Id} - K(z,e) + B\right)^{-1}$ is holomorphic in $z$ near $z_0$ and depends continuously on $e$, so the problem is reduced to proving the meromorphic invertibility of $\left(\textrm{Id} - G(z,e)\right)$ and the continuity with respect to the parameter $e$. Recall that $G(z,e)$ is meromorphic in $z$, continuous in $e$, with finite rank, so we can always represent it as
$$ 
G(z,e) = \sum_{ 1\leqslant i,j\leqslant p} a_{ij}(z,e)\,\vert\varphi_i > < \varphi_j \vert  
$$
where the coefficients $a_{ij}(z,e)$ are meromorphic in $z$, continuous in $e$ and $(\varphi_i)_{i=1}^p$ is a finite family of linearly independent vectors in $\mathcal{H}$. To solve $\big(\textrm{Id} - G(z,e)\big)v=w $ where $w$ is given, we make the ansatz $v=w+\sum_{i=1}^p b_i\varphi_i$ therefore the equation becomes
$$
\big(\textrm{Id}-G(z,e)\big)v = \big(\textrm{Id} - G(z,e)\big)\,\left(w+\sum_{i=1}^p b_i\varphi_i\right)
$$

\begin{equation*} \begin{split}
&= w+\sum_{i=1}^p b_i\varphi_i- \sum_{ 1\leqslant i,j\leqslant p,k} b_ka_{ij}(z,e)\varphi_i \left\langle \varphi_j,\varphi_k\right\rangle -\sum_{ 1\leqslant i,j\leqslant p} a_{ij}(z,e)\varphi_i \left\langle \varphi_j,w\right\rangle  
\end{split} \end{equation*}
that simplifies into the simpler relation
$$
\sum_{i=1}^p b_i\varphi_i- \sum_{ 1\leqslant i,j\leqslant p,k} b_ka_{ij}(z,u)\varphi_i \left\langle \varphi_j,\varphi_k\right\rangle =\sum_{ 1\leqslant i,j\leqslant p} a_{ij}(z,e)\varphi_i \left\langle \varphi_j,w\right\rangle.
$$
By linear algebra, the above equation can be solved on the complement of the zero locus of the polynomial
$$
\det\left(\delta_{ik}-\sum_j a_{ij}(z,e)\left\langle \varphi_j,\varphi_k\right\rangle  \right) 
$$
which depends meromorphically on $z$ and continuously on $e$. So away from the zero locus of the determinant we can meromorphically invert $\textrm{Id} - G(z,e)$ hence $\textrm{Id} -K(z,e)$ and everything depends continuously on the parameter $e$. The fact that the poles have finite rank comes from the fact that they only appear through the finite rank operator $G(z,e)$.
\end{proof}

\section{Geometric Littlewood-Paley decomposition}
\label{SubsectinoGeometricLP}

We recall from Klainerman and Rodnianski's work \cite{KR} the basics of Littlewood-Paley decomposition in a manifold setting. We use it to indicate a proof of Proposition \ref{PropCvgce} on the renormalization of ${\sf \Pi}(\xi_\varepsilon, X_\varepsilon)$ which is used in the construction of the resolvent of $H$. 

\medskip

\begin{theorem}[Klainerman-Rodnianski]
Given $\ell\in\bbN$ there exists a Schwartz function $m$ such that 
\begin{equation} \label{EqNullMoments}
\int_0^\infty t^{k_1}\partial_t^{k_2} m(t)\,dt = 0\qquad (\forall (k_1,k_2),\; k_1+k_2\leqslant \ell)
\end{equation}
and such that the self--adjoint smoothing operators
\begin{equation} \label{EqDefnPk}
P_k = \int_{0}^\infty 2^{2k} m(2^{2k}t) e^{t\Delta} \, dt  \qquad (k\in \bbN\cup\{-1\})
\end{equation}
enjoy the following properties.   \vspace{0.15cm}

\begin{enumerate}
   \item[{\sf (a)}] \textbf{Resolution of the identity.} One has $\sum_{k\geq -1} P_k = \textrm{\emph{Id}}$.   \vspace{0.1cm}
 
   \item[{\sf (b)}] \textbf{Bessel inequality.} One has
$$ 
\sum_{k\geq 0} \Vert P_kf \Vert_{L^2}\lesssim \Vert f\Vert_{L^2}.
$$  

   \item[{\sf (c)}] \textbf{Finite band property.} One has 
   $$
   \Vert \Delta P_kf\Vert_{L^p}\lesssim 2^{2k}\Vert f \Vert_{L^p},
   $$
   and
   $$
   \Vert P_kf\Vert_{L^p}\lesssim 2^{-2k}\Vert \Delta f\Vert_{L^p};
   $$
   also we have the dual estimate $\Vert P_k\nabla f \Vert_{L^2}\lesssim 2^k\Vert f\Vert_{L^2}$,   \vspace{0.1cm}
   
   \item[{\sf (d)}] \textbf{Flexibility property.} There exists a function $\tilde{m}$ satisfying \eqref{EqNullMoments} such that $\Delta P_k=2^{2k}\tilde{P}_k$ and the family $(\tilde{P}_k)_k$ is a Littlewood--Paley decomposition which might not satisfy the resolution of identity equation.
\end{enumerate}
\end{theorem}

\medskip

We quickly recall the main features of the heat calculus we shall use in the sequel. The heat calculus is a way to encode the salient features of the Euclidean heat kernel $(4\pi t)^{-\frac{d}{2}}e^{-\frac{\Vert x-y \Vert^2}{4t}}$ and of the first approximation of the heat kernel on manifolds $K_1(t,x,y) = (4\pi t)^{\frac{d}{2}}e^{-\frac{\Vert x-y \Vert^2_{g(y)}}{4t}}$, which are
\begin{itemize}
   \item the prefactor $t^{-\frac{d}{2}}$,
   \item the exponential factor, which is a smooth function of $X=\frac{x-y}{\sqrt{t}}$ and
$y$, exponentially decaying as $\Vert X\Vert\rightarrow +\infty$.
\end{itemize}
This motivates the following definition, in which the notation $C^\infty([0,+\infty)_{\frac{1}{2}})$ stands for the set of functions $f(t)$ which are smooth as functions of $\sqrt{t}$, for $t\geqslant 0$.

\medskip

\begin{definition} \label{def:heat} 
Pick a non-positive index $\gamma$. The space $\Psi^\gamma_H$ is defined to be the set of functions in $C^\infty((0,+\infty)\times S^2)$ satisfying the following axioms
\begin{itemize}
\item $A$ is smooth, if $x\neq y$ then $A(t,x,y) = O(t^\infty)$,
\item For any $p\in M$, there exists a chart $U$ containing $p$ and $\tilde{A}\in C^\infty\big([0,+\infty)_{\frac{1}{2}}\times U\times \mathbb{R}^d\big)$ such that for $(x,y)\in U^2$ one has 
$$
A(t,x,y)=t^{-\frac{d+2}{2}-\gamma}\widetilde{A}\Big(\sqrt{t},\frac{x-y}{\sqrt{t}},y\Big) 
$$
where $\widetilde{A}$ has rapid decay in the second variable
$$
\big\Vert D^\gamma_{\sqrt{t},X,y} \tilde{A} \big\Vert = O\big(\Vert X\Vert^{-\infty}\big)
$$
when $\Vert X\Vert\rightarrow +\infty$. 
\end{itemize}
\end{definition}

\medskip

The use of the heat calculus gives a familiar form to the operators $P_k$. Set 
$$
M(t) \defeq \int_0^tm(s)ds
$$ 
and use the presentation of the heat calculus in the chart from definition~\ref{def:heat} to write
\begin{eqnarray*}
\int_0^\infty  2^jm(2^jt)e^{t\Delta}(x,y)dt=2^{-j}
\int_0^\infty  M(t) 2^{k\frac{d}{2}}t^{-\frac{d}{2}}\tilde{A}\Big(2^{-k}t,x,2^{\frac{k}{2}}\frac{x-y}{\sqrt{t}}\Big) \, dt.
\end{eqnarray*}
Then for any pair of test functions $\chi_1,\chi_2$
\begin{eqnarray*}
\widehat{(P_k\chi_1)\chi_2}(\xi,\eta)&=&2^{-k}\int_{U\times\mathbb{R}^2}\chi_1(x)\chi_2(h)e^{i(\xi.x+h.\eta)}\int_0^\infty  M(t) 2^{k\frac{d}{2}}t^{-\frac{d}{2}}\tilde{A}\Big(2^{-k}t,x,2^{\frac{k}{2}}\frac{h}{\sqrt{t}}\Big) \, dtdx2h   \\
&=&2^{-k}\int_{U\times\mathbb{R}^2}\chi_1(x)\chi_2(2^{-\frac{k}{2}}h) e^{i(\xi.x+2^{-\frac{k}{2}}h.\eta)}\int_0^\infty  M(t) t^{-\frac{d}{2}}\tilde{A}\Big(2^{-k}t,x,\frac{h}{\sqrt{t}}\Big) \, dtdxdh.
\end{eqnarray*}
Using the rapid decay in the $h$ variable for all values of $t\in [0,+\infty)$, $x\in U$
$$
\sup_{x\in U} \Big\vert \tilde{A}\Big(2^{-k}t,x,\frac{h}{\sqrt{t}}\Big) \Big\vert \leqslant C_N\big(1+\vert h\vert\big)^{-N} 
$$ 
and the fact that
$\chi_1(x)\chi_2(2^{-\frac{k}{2}}h)\int_0^\infty  M(t) t^{-\frac{d}{2}}\tilde{A}\big(2^{-k}t,x,\frac{h}{\sqrt{t}}\big)\,dt$ is bounded in $C^\infty(U\times \mathbb{R}^2)$ uniformly in the parameter $k$, we have an estimate of the form
\begin{eqnarray*}
\big\vert \widehat{(P_k\chi_1)\chi_2}(\xi,\eta)\big\vert \leqslant C_N 2^{-k} \left(1+\vert\xi\vert+2^{-\frac{k}{2}}\vert\eta \vert \right)^{-N}
\end{eqnarray*}
In position space, in the local chart $U\times U$ from definition~\ref{def:heat}, the estimate reads
\begin{equation} \label{EqRepresentationLPProjector}
P_k(x,y) = 2^{-k} 2^{k\frac{d}{2}} K_k\big(x,2^{\frac{k}{2}}(x-y)\big),
\end{equation}
where the $(K_k)_k$ form a bounded family of smooth functions in $C^\infty\big(U\times \big\{\vert h\vert\leqslant 1\big\}\big)$. 

\medskip

Let $P$ and $\widetilde{P}$ be a family of geometric Littlewood-Paley projectors built from functions $m$ and $\widetilde m$ that vanish at $t=0$. It will be convenient in the proof of Proposition \ref{PropCvgce} to control the kernel $\sum_{i,j\geqslant 0} \big((\Delta^\alpha P_i)P_j\big)(x,y)$ in terms of $\alpha$. We know from p.140 of \cite{KR} that we have the exact identity
$$
\big({\widetilde P}_iP_j\big)(x,y) = - 2^{-2\vert i-j\vert}\int_0^\infty\int_0^\infty \int_0^1  e^{(t_1+st_2)\Delta}{\widetilde m}_i(t_1) t_2m_j(t_2)\,dsdt_1dt_2.
$$
Using the structure of the heat kernel which follows from the heat calculus we may write in local coordinate chart $x\in U, h\in \mathbb{R}^2$
\begin{eqnarray*}
\big({\widetilde P}_iP_j\big)(x,x+h) \eqdef 2^{-2\vert i-j\vert}K_{ij}(x,h)
\end{eqnarray*}
where 
\begin{eqnarray}\label{e:estimatekernel}
\sup_{x\in U} \Big\vert \partial_h^\beta\partial_x^\alpha K_{ij}(x,h) \Big\vert \leqslant C_{\alpha,\beta} \, 2^{-i}2^{i\frac{d}{2}} 2^{i\frac{\vert \beta\vert}{2}}
\end{eqnarray}
uniformly in $(i,j)$. These are the seminorms for the topology of distributions whose wavefront set is concentrated on the conormal bundle of the diagonal.

Now in~\cite{KR} we also find that $\Delta^\alpha P_i\widetilde{P}_j=2^{2i\alpha}Q_iP_j$, where $(Q_i)_i$ is an admissible family of Littlewood-Paley projectors. We deduce from this observation an estimate of the form
\begin{eqnarray*}
(\Delta^\alpha P_i) P_j(x,x+h) = 2^{2i\alpha}2^{-2\vert i-j\vert}K_{ij}(x,h),
\end{eqnarray*} 
where the kernel $K_{ij}$ satisfies the same estimate \ref{e:estimatekernel}. This is all we need to prove the following technical lemma.

\medskip

\begin{lemma} \label{LemEstimatesRenormalization}
Let the Littlewood-Paley projectors $(P_i)_i$ be constructed from a function $m$ that vanishes at $t=0$. Fix $k\geq 1$ and $(\alpha_1,\dots,\alpha_k)\in\bbZ^k$. The series of Schwartz kernels
\begin{eqnarray} \label{EqSumTechnicalLemma}
\sum_{i_1,\dots,i_k,j_1,\dots,j_k}\sum_{\vert i_1-i_2\vert \leq 1,\dots,\vert i_1-i_k\vert\leqslant 1} \big((\Delta^{\alpha_1}P_{i_1})P_{j_1}\big)(x,y)\dots \big((\Delta^{\alpha_k}P_{i_k})P_{j_k}\big)(x,y) 
\end{eqnarray}
converges absolutely in the space of pseudodifferential kernels of order $2(\alpha_1+\dots+\alpha_k)+(k-1)\frac{d}{2}$.
\end{lemma}

\medskip

\begin{proof}
Using the above discussion we may rewrite
\begin{equation*} \begin{split}
\big((\Delta^{\alpha_1}P_{i_1})P_{j_1}\big)&(x,y)\dots \big((\Delta^{\alpha_k}P_{i_k})P_{j_k}\big)(x,y)   \\
&= 2^{2(i_1\alpha_1+\dots+i_k\alpha_k)}2^{-2(\vert i_1-j_1\vert+\dots+\vert i_k-j_k \vert)}K_{i_1j_1}(x,y)\dots K_{i_kj_k}(x,y)
\end{split} \end{equation*}
where the smooth functions $ K_{i_n j_n}(x,y)$ satisfy the estimate \eqref{e:estimatekernel}. So one has for all tuples $(i_1,\dots,i_k,j_1,\dots,j_k)$ such that  $\vert i_1-i_2\vert \leq 1,\dots,\vert i_1-i_k\vert\leqslant 1$ an estimate of the form

\begin{equation*} \begin{split}
\Big\vert \partial_h^b \partial_x^a K_{i_1,j_1}(x,x+h)\dots K_{i_k,j_k}&(x,x+h)\Big\vert   \\
&\leqslant  C_{ab} \, 2^{-(i_1+\dots+i_k)} 2^{(i_1+\dots+i_k)\frac{d}{2}}   2^{2 \inf(i_1,\dots,i_k,j_1,\dots,j_k) \frac{\vert b\vert}{2}}
\end{split} \end{equation*}
where the constant $C_{ab}$ does not depend on the indices $(i_1,\dots,i_k,j_1,\dots,j_k)$. This estimate ensures that the sum \eqref{EqSumTechnicalLemma} converges in the space of conormal distributions of order $2(\alpha_1+\dots+\alpha_k)+(k-1)\frac{d}{2}$.
\end{proof}

\medskip

We give here the proof of Proposition \ref{PropCvgce} performing the Wick renormalization of the resonant term ${\sf \Pi}(\xi_\varepsilon , X_\varepsilon)$.

\medskip

\begin{proof}
\textit{Step 1 -- Singular part.} Since the two paraproduct terms in the decomposition of the product $\xi_\varepsilon X_{\varepsilon}$ converge as $\varepsilon$ goes to $0$ the quantities $\bbE\big[{\sf \Pi}(\xi_\varepsilon,X_{\varepsilon})\big]$ and $\bbE\big[\xi_\varepsilon X_{\varepsilon}\big]$ differ by a convergent quantity. Use now the Markov property of the heat operator and the definition of white noise to see that 
$$
\mathbb{E}\left[\left(\left(-\Delta +z_0\right)^{-1}\xi_\varepsilon\right)(x)\xi_\varepsilon (x)\right]=  \big( e^{2\varepsilon\Delta}\left(-\Delta +z_0\right)^{-1}\big)(x,x).
$$ 
 So the singular part of the above expectation
$$
\mathbb{E}\left[\left(\left(-\Delta +z_0\right)^{-1}\xi_\varepsilon\right)(x)\xi_\varepsilon (x)\right]
$$ 
comes from the term $\left(e^{2\varepsilon\Delta} \left(-\Delta +z_0\right)^{-1}\right)(x,x)$.

 An immediate computation yields
$$
e^{2\varepsilon\Delta}\left(-\Delta +z_0\right)^{-1} = \int_{2\varepsilon}^1 e^{(-z_0-2r)s} e^{s\Delta}(\textrm{Id} - \pi_0) ds+\int_1^\infty e^{s\Delta}(\textrm{Id} - \pi_0)e^{-z_0s} ds
$$
where $\pi_0$ is the orthogonal projector on the subspace of constant functions. Recall that $z_0$ is very positive so the integral over $[1,\infty)$ converges absolutely and defines a smoothing operator; it does not contribute to the singular part of $\big(e^{2\varepsilon\Delta}\left(-\Delta + z_0\right)^{-1}\big)(x,x)$ when $\varepsilon$ goes to $0$. Now using the asymptotic expansion of the heat kernel yields the identity 
$$
\big(e^{s\Delta}(\textrm{Id}-\pi_0)\big)(x,x) = \frac{1}{4\pi s}+O(1),
$$ 
with an error term $O(1)$ bounded in $s$ and smooth in the $x$ variable. It follows that
\begin{eqnarray*}
\big(e^{2\varepsilon\Delta}\left(-\Delta + z_0\right)^{-1}\big)(x,x)=  \int_{2\varepsilon}^1e^{(-z_0-2\varepsilon)s} \frac{1}{4\pi s} ds+\mathcal{O}(1) = \frac{\vert \log(\varepsilon)\vert}{4\pi} + O(1).
\end{eqnarray*}
We see here that the singular part of $\bbE\big[\xi_\varepsilon X_{\varepsilon}\big]$ does not depend on the point $x$.

\textit{Step 2 -- Stochastic estimates.} 
For the Kolmogorov type estimates, we refer to~\cite{BDFT} where these estimates are done in detail and the hypercontractivity is applied carefully. The difficulty in the curved case is that we lost the stationarity of the law of the process.
\end{proof}

\bigskip

\noindent \textcolor{gray}{$\bullet$} I. Bailleul -- Univ Brest, CNRS UMR 6205, Laboratoire de Mathematiques de Bretagne Atlantique, France.   \\
\noindent {\it E-mail}: ismael.bailleul@univ-brest.fr 

\smallskip

\noindent \textcolor{gray}{$\bullet$} V. Dang -- Sorbonne Université -- Université de Paris, CNRS, UMR 7586, Paris, France. 
Institut Universitaire de France.    \\
{\it E-mail}: nguyen-viet.dang@imj-prg.fr

\smallskip

\noindent \textcolor{gray}{$\bullet$} A. Mouzard -- Modal’X - UMR CNRS 9023, Université Paris Nanterre, 92000 Nanterre, France   \\
{\it E-mail}: antoine.mouzard@math.cnrs.fr

\end{document}